\documentclass[12pt]{amsart}
\usepackage{amssymb}
\begin{document}
\theoremstyle{plain}
\newtheorem{thm}{Theorem}[section]
\newtheorem{theorem}[thm]{Theorem}
\newtheorem{lemma}[thm]{Lemma}
\newtheorem{corollary}[thm]{Corollary}
\newtheorem{proposition}[thm]{Proposition}
\newtheorem{addendum}[thm]{Addendum}
\newtheorem{variant}[thm]{Variant}
\theoremstyle{definition}
\newtheorem{notations}[thm]{Notations}
\newtheorem{question}[thm]{Question}
\newtheorem{problem}[thm]{Problem}
\newtheorem{remark}[thm]{Remark}
\newtheorem{remarks}[thm]{Remarks}
\newtheorem{definition}[thm]{Definition}
\newtheorem{claim}[thm]{Claim}
\newtheorem{assumption}[thm]{Assumption}
\newtheorem{assumptions}[thm]{Assumptions}
\newtheorem{properties}[thm]{Properties}
\newtheorem{example}[thm]{Example}
\numberwithin{equation}{thm}
\catcode`\@=11
\def\opn#1#2{\def#1{\mathop{\kern0pt\fam0#2}\nolimits}}
\def\bold#1{{\bf #1}}%
\def\underrightarrow{\mathpalette\underrightarrow@}
\def\underrightarrow@#1#2{\vtop{\ialign{$##$\cr
 \hfil#1#2\hfil\cr\noalign{\nointerlineskip}%
 #1{-}\mkern-6mu\cleaders\hbox{$#1\mkern-2mu{-}\mkern-2mu$}\hfill
 \mkern-6mu{\to}\cr}}}
\let\underarrow\underrightarrow
\def\underleftarrow{\mathpalette\underleftarrow@}
\def\underleftarrow@#1#2{\vtop{\ialign{$##$\cr
 \hfil#1#2\hfil\cr\noalign{\nointerlineskip}#1{\leftarrow}\mkern-6mu
 \cleaders\hbox{$#1\mkern-2mu{-}\mkern-2mu$}\hfill
 \mkern-6mu{-}\cr}}}
\let\amp@rs@nd@\relax
\newdimen\ex@
\ex@.2326ex
\newdimen\bigaw@
\newdimen\minaw@
\minaw@16.08739\ex@
\newdimen\minCDaw@
\minCDaw@2.5pc
\newif\ifCD@
\def\minCDarrowwidth#1{\minCDaw@#1}
\newenvironment{CD}{\@CD}{\@endCD}
\def\@CD{\def\A##1A##2A{\llap{$\vcenter{\hbox
 {$\scriptstyle##1$}}$}\Big\uparrow\rlap{$\vcenter{\hbox{%
$\scriptstyle##2$}}$}&&}%
\def\V##1V##2V{\llap{$\vcenter{\hbox
 {$\scriptstyle##1$}}$}\Big\downarrow\rlap{$\vcenter{\hbox{%
$\scriptstyle##2$}}$}&&}%
\def\={&\hskip.5em\mathrel
 {\vbox{\hrule width\minCDaw@\vskip3\ex@\hrule width
 \minCDaw@}}\hskip.5em&}%
\def\verteq{\Big\Vert&&}%
\def\noarr{&&}%
\def\vspace##1{\noalign{\vskip##1\relax}}\relax\let\amp@rs@nd@&\iffalse}\fi
 \CD@true\vcenter\bgroup\relax\let\\=\cr\iffalse}\fi\tabskip\z@skip\baselineskip20\ex@
 \lineskip3\ex@\lineskiplimit3\ex@\halign\bgroup
 &\hfill$\m@th##$\hfill\cr}
\def\@endCD{\cr\egroup\egroup}
\def\>#1>#2>{\amp@rs@nd@\setbox\z@\hbox{$\scriptstyle
 \;{#1}\;\;$}\setbox\@ne\hbox{$\scriptstyle\;{#2}\;\;$}\setbox\tw@
 \hbox{$#2$}\ifCD@
 \global\bigaw@\minCDaw@\else\global\bigaw@\minaw@\fi
 \ifdim\wd\z@>\bigaw@\global\bigaw@\wd\z@\fi
 \ifdim\wd\@ne>\bigaw@\global\bigaw@\wd\@ne\fi
 \ifCD@\hskip.5em\fi
 \ifdim\wd\tw@>\z@
 \mathrel{\mathop{\hbox to\bigaw@{\rightarrowfill}}\limits^{#1}_{#2}}\else
 \mathrel{\mathop{\hbox to\bigaw@{\rightarrowfill}}\limits^{#1}}\fi
 \ifCD@\hskip.5em\fi\amp@rs@nd@}
\def\<#1<#2<{\amp@rs@nd@\setbox\z@\hbox{$\scriptstyle
 \;\;{#1}\;$}\setbox\@ne\hbox{$\scriptstyle\;\;{#2}\;$}\setbox\tw@
 \hbox{$#2$}\ifCD@
 \global\bigaw@\minCDaw@\else\global\bigaw@\minaw@\fi
 \ifdim\wd\z@>\bigaw@\global\bigaw@\wd\z@\fi
 \ifdim\wd\@ne>\bigaw@\global\bigaw@\wd\@ne\fi
 \ifCD@\hskip.5em\fi
 \ifdim\wd\tw@>\z@
 \mathrel{\mathop{\hbox to\bigaw@{\leftarrowfill}}\limits^{#1}_{#2}}\else
 \mathrel{\mathop{\hbox to\bigaw@{\leftarrowfill}}\limits^{#1}}\fi
 \ifCD@\hskip.5em\fi\amp@rs@nd@}
\newenvironment{CDS}{\@CDS}{\@endCDS}
\def\@CDS{\def\A##1A##2A{\llap{$\vcenter{\hbox
 {$\scriptstyle##1$}}$}\Big\uparrow\rlap{$\vcenter{\hbox{%
$\scriptstyle##2$}}$}&}%
\def\V##1V##2V{\llap{$\vcenter{\hbox
 {$\scriptstyle##1$}}$}\Big\downarrow\rlap{$\vcenter{\hbox{%
$\scriptstyle##2$}}$}&}%
\def\={&\hskip.5em\mathrel
 {\vbox{\hrule width\minCDaw@\vskip3\ex@\hrule width
 \minCDaw@}}\hskip.5em&}
\def\verteq{\Big\Vert&}
\def\novarr{&}
\def\noharr{&&}
\def\SE##1E##2E{\slantedarrow(0,18)(4,-3){##1}{##2}&}
\def\SW##1W##2W{\slantedarrow(24,18)(-4,-3){##1}{##2}&}
\def\NE##1E##2E{\slantedarrow(0,0)(4,3){##1}{##2}&}
\def\NW##1W##2W{\slantedarrow(24,0)(-4,3){##1}{##2}&}
\def\slantedarrow(##1)(##2)##3##4{%
\thinlines\unitlength1pt\lower 6.5pt\hbox{\begin{picture}(24,18)%
\put(##1){\vector(##2){24}}%
\put(0,8){$\scriptstyle##3$}%
\put(20,8){$\scriptstyle##4$}%
\end{picture}}}
\def\vspace##1{\noalign{\vskip##1\relax}}\relax\let\amp@rs@nd@&\iffalse}\fi
 \CD@true\vcenter\bgroup\relax\let\\=\cr\iffalse}\fi\tabskip\z@skip\baselineskip20\ex@
 \lineskip3\ex@\lineskiplimit3\ex@\halign\bgroup
 &\hfill$\m@th##$\hfill\cr}
\def\@endCDS{\cr\egroup\egroup}
\newdimen\TriCDarrw@
\newif\ifTriV@
\newenvironment{TriCDV}{\@TriCDV}{\@endTriCD}
\newenvironment{TriCDA}{\@TriCDA}{\@endTriCD}
\def\@TriCDV{\TriV@true\def\TriCDpos@{6}\@TriCD}
\def\@TriCDA{\TriV@false\def\TriCDpos@{10}\@TriCD}
\def\@TriCD#1#2#3#4#5#6{%
\setbox0\hbox{$\ifTriV@#6\else#1\fi$}
\TriCDarrw@=\wd0 \advance\TriCDarrw@ 24pt
\advance\TriCDarrw@ -1em
\def\SE##1E##2E{\slantedarrow(0,18)(2,-3){##1}{##2}&}
\def\SW##1W##2W{\slantedarrow(12,18)(-2,-3){##1}{##2}&}
\def\NE##1E##2E{\slantedarrow(0,0)(2,3){##1}{##2}&}
\def\NW##1W##2W{\slantedarrow(12,0)(-2,3){##1}{##2}&}
\def\slantedarrow(##1)(##2)##3##4{\thinlines\unitlength1pt
\lower 6.5pt\hbox{\begin{picture}(12,18)%
\put(##1){\vector(##2){12}}%
\put(-4,\TriCDpos@){$\scriptstyle##3$}%
\put(12,\TriCDpos@){$\scriptstyle##4$}%
\end{picture}}}
\def\={\mathrel {\vbox{\hrule
   width\TriCDarrw@\vskip3\ex@\hrule width
   \TriCDarrw@}}}
\def\>##1>>{\setbox\z@\hbox{$\scriptstyle
 \;{##1}\;\;$}\global\bigaw@\TriCDarrw@
 \ifdim\wd\z@>\bigaw@\global\bigaw@\wd\z@\fi
 \hskip.5em
 \mathrel{\mathop{\hbox to \TriCDarrw@
{\rightarrowfill}}\limits^{##1}}
 \hskip.5em}
\def\<##1<<{\setbox\z@\hbox{$\scriptstyle
 \;{##1}\;\;$}\global\bigaw@\TriCDarrw@
 \ifdim\wd\z@>\bigaw@\global\bigaw@\wd\z@\fi
 \mathrel{\mathop{\hbox to\bigaw@{\leftarrowfill}}\limits^{##1}}
 }
 \CD@true\vcenter\bgroup\relax\let\\=\cr\iffalse}\fi
 \tabskip\z@skip\baselineskip20\ex@
 \lineskip3\ex@\lineskiplimit3\ex@
 \ifTriV@
 \halign\bgroup
 &\hfill$\m@th##$\hfill\cr
#1&\multispan3\hfill$#2$\hfill&#3\\
&#4&#5\\
&&#6\cr\egroup%
\else
 \halign\bgroup
 &\hfill$\m@th##$\hfill\cr
&&#1\\%
&#2&#3\\
#4&\multispan3\hfill$#5$\hfill&#6\cr\egroup
\fi}
\def\@endTriCD{\egroup}
\newcommand{\sA}{{\mathcal A}}
\newcommand{\sB}{{\mathcal B}}
\newcommand{\sC}{{\mathcal C}}
\newcommand{\sD}{{\mathcal D}}
\newcommand{\sE}{{\mathcal E}}
\newcommand{\sF}{{\mathcal F}}
\newcommand{\sG}{{\mathcal G}}
\newcommand{\sH}{{\mathcal H}}
\newcommand{\sI}{{\mathcal I}}
\newcommand{\sJ}{{\mathcal J}}
\newcommand{\sK}{{\mathcal K}}
\newcommand{\sL}{{\mathcal L}}
\newcommand{\sM}{{\mathcal M}}
\newcommand{\sN}{{\mathcal N}}
\newcommand{\sO}{{\mathcal O}}
\newcommand{\sP}{{\mathcal P}}
\newcommand{\sQ}{{\mathcal Q}}
\newcommand{\sR}{{\mathcal R}}
\newcommand{\sS}{{\mathcal S}}
\newcommand{\sT}{{\mathcal T}}
\newcommand{\sU}{{\mathcal U}}
\newcommand{\sV}{{\mathcal V}}
\newcommand{\sW}{{\mathcal W}}
\newcommand{\sX}{{\mathcal X}}
\newcommand{\sY}{{\mathcal Y}}
\newcommand{\sZ}{{\mathcal Z}}
\newcommand{\A}{{\mathbb A}}
\newcommand{\B}{{\mathbb B}}
\newcommand{\C}{{\mathbb C}}
\newcommand{\D}{{\mathbb D}}
\newcommand{\E}{{\mathbb E}}
\newcommand{\F}{{\mathbb F}}
\newcommand{\G}{{\mathbb G}}
\newcommand{\HH}{{\mathbb H}}
\newcommand{\I}{{\mathbb I}}
\newcommand{\J}{{\mathbb J}}
\newcommand{\M}{{\mathbb M}}
\newcommand{\N}{{\mathbb N}}
\renewcommand{\P}{{\mathbb P}}
\newcommand{\Q}{{\mathbb Q}}
\newcommand{\R}{{\mathbb R}}
\newcommand{\T}{{\mathbb T}}
\newcommand{\U}{{\mathbb U}}
\newcommand{\V}{{\mathbb V}}
\newcommand{\W}{{\mathbb W}}
\newcommand{\X}{{\mathbb X}}
\newcommand{\Y}{{\mathbb Y}}
\newcommand{\Z}{{\mathbb Z}}
\title[Base spaces of non-isotrivial families]{Base spaces of non-isotrivial
families of smooth minimal models}
\author[Eckart Viehweg]{Eckart Viehweg}
\address{Universit\"at Essen, FB6 Mathematik, 45117 Essen, Germany}
\email{ viehweg@uni-essen.de}
\thanks{This work has been supported by the ``DFG-Forschergruppe
Arithmetik und Geometrie'' and the ``DFG-Schwerpunktprogramm
Globale Methoden in der Komplexen Geometrie''. The second named author
is supported by a grant from the Research
Grants Council of the Hong Kong
Special Administrative Region, China
(Project No. CUHK 4239/01P).}
\author[Kang Zuo]{Kang Zuo}
\address{The Chinese University of Hong Kong, Department of Mathematics,
Shatin, Hong Kong}
\email{kzuo@math.cuhk.edu.hk}
\dedicatory{F\"ur Hans Grauert, mit tiefer Bewunderung.}
\maketitle
Given a polynomial $h$ of degree $n$ let $\sM_h$ be the moduli functor of
canonically polarized complex manifolds with Hilbert polynomial
$h$. By \cite{Vie} there exist a quasi-projective scheme $M_h$
together with a natural transformation
$$\Psi:\sM_h \to {\rm Hom}(\underline{ \ \ },M_h)$$
such that $M_h$ is a coarse moduli scheme for $\sM_h$.
For a complex quasi-projective manifold $U$ we will
say that a morphism $\varphi: U \to M_h$ factors through the
moduli stack, or that $\varphi$ is induced by a family,
if $\varphi$ lies in the image of $\Psi(U)$,
hence if $\varphi=\Psi(f:V\to U)$.

Let $Y$ be a projective non-singular compactification of $U$
such that $S=Y\setminus U$ is a normal crossing divisor,
and assume that the morphism $\varphi:U\to M_h$, induced by a
family, is generically finite. For moduli of curves of genus
$g\geq 2$, i.e. for $h(t)=(2t-1)(g-1)$, it is easy to show, that
the existence of $\varphi$ forces $\Omega^1_{Y}(\log S)$ to be
big (see \ref{0.1} for the definition), hence that
$S^m(\Omega^1_{Y}(\log S))$ contains an ample subsheaf of full
rank for some $m>0$. In particular, $U$ should not be an abelian variety
or $\C^*$. By \cite{Lu} the bigness of $\Omega^1_{Y}(\log S)$
implies even the Brody hyperbolicity of $U$. As we will see,
there are other restrictions on $U$, as those formulated below
in \ref{cor1}, \ref{cor2}, and \ref{cor5}.

In the higher dimensional case, i.e. if $\deg(h) >1$,
L. Migliorini, S. Kov\'acs, E. Bedulev and the authors studied
in \cite{Mig}, \cite{Kov1}, \cite{Kov2}, \cite{Kov3},
\cite{B-V}, \cite{V-Z} \cite{V-Z2} geometric properties of
manifolds $U$ mapping non-trivially to the moduli stack.
Again, $U$ can not be $\C^*$, nor an abelian variety, and
more generally it must be Brody hyperbolic, if $\varphi$ is
quasi-finite.

In general the sheaf $\Omega^1_{X/Y}(\log S)$
fails to be big (see example \ref{5.3}).
Nevertheless, building up on the methods introduced in
\cite{V-Z} and \cite{V-Z2} we will show that
for $m$ sufficiently large the sheaf $S^m(\Omega_{Y}^1(\log S))$
has enough global sections (see section \ref{0} for the precise
statement), to exclude the existence of a generically finite
morphism $\varphi:U \to M_h$, or even of a non-trivial morphism,
for certain manifolds $U$.

\begin{theorem}[see \ref{4.2}, \ref{4.3} and \ref{7.2}]\label{cor1} Assume
that $U$ satisfies one of the following conditions
\begin{enumerate}
\item[a)] $U$ has a smooth projective compactification $Y$
with $S=Y\setminus U$ a normal crossing divisor and with
$T_Y(-\log S)$ weakly positive.
\item[b)] Let $H_1+ \cdots + H_\ell$ be a reduced normal crossing
divisor in $\P^N$, and $\ell< \frac{N}{2}$. For $0\leq r\leq l$ define
\begin{gather*}
H = \bigcap_{j=r+1}^\ell H_j, \ \ \ S_i = H_i|_H, \ \ \
S = \sum_{i=1}^r S_i,
\end{gather*}
and assume $U= H \setminus S$.
\item[c)] $U=\P^N\setminus S$ for a reduced normal crossing divisor
$S=S_1+ \cdots + S_\ell$ in $\P^N$, with $\ell< N.$
\end{enumerate}
Then a morphism $U \to M_h$, induced by a family, must be
trivial.
\end{theorem}

In a) the sheaf $T_Y(-\log S)$ denotes the dual of the sheaf
of one forms with logarithmic poles along $S$. The definition of
``weakly positive'' will be recalled in \ref{0.1}.
Part a) of \ref{cor1}, for $S=\emptyset$, has been shown by S.
Kov\'acs in \cite{Kov3}.

Considering $r=0$ in \ref{cor1}, b), one finds that
smooth complete intersections $U$ in $\P^N$ of
codimension $\ell< \frac{N}{2}$ do not allow a non-trivial
morphism $U \to M_h$, induced by a family.
In b) the intersection with an empty index set is supposed to be
$H=\P^N$. So for $\ell < \frac{N}{2}$ part c) follows from b).

In general, \ref{cor1}, c), will follow from the slightly stronger
statement in the second part of the next theorem. In fact, if one
chooses general linear hyperplanes $D_0, \cdots , D_N$ then
$D_0+D_1+\cdots +D_N+S$ remains a normal crossing divisor. and all
morphism
$$
\P^N\setminus(D_0+D_1+\cdots +D_N+S) \>>> M_h,
$$
induced by a family, must be trivial.

\begin{theorem}[see \ref{7.2}]\label{cor5}
\begin{enumerate}
\item[a)] Assume
that $U$ is the complement of a normal crossing divisor $S$
with strictly less than $N$ components in an $N$-dimensional
abelian variety. Then there exists no generically finite
morphism $U \to M_h$, induced by a family.
\item[b)] For $Y=\P^{\nu_1}\times \cdots \times
\P^{\nu_k}$ let
$$
D^{(\nu_i)}=D_0^{(\nu_i)}+\cdots +D_{\nu_i}^{(\nu_i)}
$$
be coordinate axes in $\P^{\nu_i}$ and
$$
D=\bigoplus_{i=1}^kD^{(\nu_i)}.
$$
Assume that $S=S_1+\cdots S_\ell$ is a divisor, such that $D+S$
is a reduced normal crossing divisor, and $\ell < \dim(Y)$.
Then there exists no morphism $\varphi:U=Y\setminus (D+S) \to M_h$
with
$$
\dim(\varphi(U)) > {\rm Max}\{\dim(Y)-\nu_i; \ i=1,\ldots ,k\}.
$$
\end{enumerate}
\end{theorem}

We do not know whether the bound $\ell < \dim(Y)$ in
\ref{cor5} is really needed. If the infinitesimal Torelli
theorem holds true for the general fibre, hence if the family
$V \to U$ induces a generically finite map to a period domain,
then the fundamental group of $U$ should not be abelian. In particular
$U$ can not be the complement of a normal crossing
divisor in $\P^N$.

In section \ref{5} we will prove different properties of $U$ in
case there exists a quasi-finite morphism $\varphi:U \to M_h$.
Those properties will be related to the rigidity of generic
curves in moduli stacks.

\begin{theorem}[see \ref{5.4} and \ref{5.7}]\label{cor2}
Let $U$ be a quasi-projective variety and let $\varphi:U \to
M_h$ be a quasi-finite morphism, induced by a family. Then
\begin{enumerate}
\item[a)] $U$ can not be isomorphic to the product of more than
$n=\deg(h)$ varieties of positive dimension.
\item[b)] ${\rm Aut}(U)$ is finite.
\end{enumerate}
\end{theorem}

Although we do not need it in its full strength, we could not
resist to include a proof of the finiteness theorem \ref{5.2}, saying that
for a projective curve $C$, for an open sub curve $C_0$, and for
a projective compactification $\bar{U}$ of $U$, the
morphisms $\pi:C \to \bar{U}$ with $\pi(C_0)\subset U$ are
parameterized by a scheme of finite type.\\

We call $f:V\to U$ a (flat or smooth) family of projective
varieties, if $f$ is projective (flat or smooth) and all fibres
connected. For a flat family, an invertible sheaf $\sL$ on $V$
will be called $f$-semi-ample, or relatively semi-ample over
$U$, if for some $\nu >0$ the evaluation of sections
$f^*f_*\sL^\nu \to \sL^\nu$ is surjective. The notion
$f$-ampleness will be used if in addition for $\nu\gg 0$ the
induced $U$-morphism $V \to \P(f_*\sL^\nu)$ is an embedding, or
equivalently, if the restriction of $\sL$ to all the fibres is
ample.

For families over a higher dimensional base $f:V\to U$, the
non-isotriviality will be measured by an invariant, introduced in
\cite{Vie1}. We define ${\rm Var}(f)$ to be the smallest
integer $\eta$ for which there exists a finitely generated
subfield $K$ of $\overline{\C(U)}$ of transcendence degree
$\eta$ over $\C$, a variety $F'$ defined over $K$, and a
birational equivalence
$$
V\times_U {\rm Spec}(\overline{\C(U)}) \sim F'\times_{{\rm Spec}(K)}
{\rm Spec}(\overline{\C(U)}).
$$
We will call $f$ isotrivial, in case that ${\rm Var}(f)=0$.
If $(f:V\to U)\in \sM_h(U)$ induces the morphism
$\varphi:U \to M_h$, then ${\rm Var}(f)=\dim(\varphi(U))$.\\

Most of the results in this article carry over to families
$V\to U$ with $\omega_{V/U}$ semi-ample. The first result
without requiring local Torelli theorems, saying that there are no
non-isotrivial families of elliptic surfaces over $\C^*$ or over
elliptic curves, has been shown by K. Oguiso and the first named
author \cite{O-V}. It was later extended to all families of higher
dimensional minimal models in \cite{V-Z}.

\begin{variant}\label{cor3}
Let $U$ be a quasi-projective manifold as in \ref{cor1}
or in \ref{cor5}. Then there exists no smooth family $f:V\to U$ with
$\omega_{V/U}$ $f$-semi-ample and with ${\rm Var}(f)=\dim(U)$.
\end{variant}

\begin{variant}\label{cor4}
For $U$ a quasi-projective manifold let $f:V\to U$ be a smooth
family with $\omega_{V/U}$ $f$-semi-ample and with ${\rm
Var}(f)=\dim(U)$. Then the conclusion a) and b) in \ref{cor2}
hold true.
\end{variant}

All the results mentioned will be corollaries of
theorem \ref{0.3}, formulated in the first section.
It is closely related to some conjectures and open problems
on differential forms on moduli stacks, explained in \ref{0.2}.
The proof of \ref{0.3}, which covers sections \ref{1}, \ref{2}, and
\ref{3}, turns out to be quite complicated, and we will try
to give an outline at the end of the first section.

The methods are close in spirit to the ones used in
\cite{V-Z} for $Y$ a curve, replacing \cite{V-Z}, Proposition 1.3, by
\cite{Zuo}, Theorem 0.1, and using some of the tools developed
in \cite{V-Z2}. So the first three and a half sections do hardly
contain any new ideas. They are needed nevertheless to adapt
methods and notations to the situation studied here, and
hopefully they can serve as a reference for methods needed to
study positivity problems over higher dimensional
bases. The reader who just wants to get some idea on the
geometry of moduli stacks should skip sections \ref{1}, \ref{2}
and \ref{3} in a first reading and start with sects \ref{0},
\ref{4}, \ref{5} and \ref{7}.\\

This article benefited from discussions between the first named author
and S. Kov\'acs. In particular we thank him for informing us
about his results. We thank the referee for hints, how to
improve the presentation of the results.

A first version of this paper was written during a visit of the
first named author to the Institute of Mathematical Science and
the Department of Mathematics at the Chinese University of Hong
Kong. The final version, including section 7, was finished during
a visit of the second named author to the Department of
Mathematics at the University of Essen. We both would like to all
the members of the host institutes for their hospitality.

\section{Differentialforms on moduli stacks}
\label{0}

Our motivation and starting point are conjectures and questions
on the sheaf of differential forms on moduli-stacks. Before
formulating the technical main result and related conjectures and
questions, let us recall some definitions.

\begin{definition} \label{0.1}
Let $\sF$ be a torsion free coherent sheaf on a quasi-projective
normal variety $Y$ and let $\sH$ be an ample invertible sheaf.
\begin{enumerate}
\item[a)] $\sF$ is generically generated if the natural morphism
$$
H^0 (Y, \sF) \otimes \sO_Y \>>> \sF
$$
is surjective over some open dense subset $U_0$ of $Y$. If one wants to
specify $U_0$ one says that $\sF$ is globally generated over $U_0$.
\item[b)] $\sF$ is weakly positive if there exists some
dense open subset $U_0$ of $Y$ with $\sF|_{U_0}$ locally free, and
if for all $\alpha > 0$ there exists some $\beta > 0$ such that
$$
S^{\alpha \cdot \beta} (\sF) \otimes \sH^{\beta}
$$
is globally generated over $U_0$. We will also say that $\sF$ is
weakly positive over $U_0$, in this case.
\item[c)] $\sF$ is big if there exists some open dense subset
$U_0$ in $Y$ and some $\mu > 0$ such that
$$
S^{\mu} (\sF) \otimes \sH^{-1}
$$
is weakly positive over $U_0$. Underlining the role of $U_0$ we will
also call $\sF$ ample with respect to $U_0$.
\end{enumerate}
\end{definition}
Here, as in \cite{Vie1} and \cite{V-Z2}, we use the following
convention: If $\sF$ is a coherent torsion free sheaf on a quasi-projective
normal variety $Y$, we consider the largest open subscheme $i: Y_1 \to
Y$ with $i^* \sF$ locally free. For
$$
\Phi = S^{\mu}, \ \ \
\Phi =\bigotimes^\mu \mbox{ \ \ \ or \ \ \ }\Phi = \det
$$
we define
$$
\Phi (\sF) = i_* \Phi (i^* \sF).
$$
Let us recall two simple properties of sheaves which are ample
with respect to open sets, or generically generated. A more
complete list of such properties can be found in \cite{Vie}, \S
2. First of all the ampleness property can be expressed in a different way
(see \cite{V-Z2}, 3.2, for example).
\begin{lemma} \label{2.1} Let $\sH$ be an ample invertible sheaf,
and $\sF$ a coherent torsion free sheaf on $Y$, whose restriction
to some open dense subset $U_0\subset Y$ is locally free. Then $\sF$ is
ample with respect to $U_0$ if and only if for some $\eta > 0$ there exists
a morphism
$$
\bigoplus \sH \>>> S^{\eta} (\sF),
$$
surjective over $U_0$.
\end{lemma}
We will also need the following well known property of generically
generated sheaves.
\begin{lemma}\label{2.2}
Let $\psi:Y'\to Y$ be a finite morphism and let $\sF$ be a coherent
torsion free sheaf on $Y$ such that $\psi^*\sF$ is generically generated.
Then for some $\beta>0$, the sheaf $S^\beta(\sF)$ is generically
generated.
\end{lemma}
\begin{proof}
We may assume that $\sF$ is locally free, and replacing $Y'$ by
some covering, that $Y'$ is a Galois cover of $Y$ with Galois
group $G$. Let
$$
\pi:\P=\P(\sF)\>>> Y\mbox{ \ \  and \ \
}\pi':\P'=\P(\psi^*\sF)\>>> Y'
$$
be the projective bundles. The induced covering $\psi':\P'\to \P$ is
again Galois. By assumption, for some $U_0\subset Y$ the sheaf
$\sO_{\P'}(1)$ is generated by global sections over
${\psi'}^{-1}\pi^{-1}(U_0)$. Hence for $g=\#G$ the sheaf
$\sO_{\P'}(g)={\psi'}^*\sO_\P(g)$ is generated over ${\psi'}^{-1}\pi^{-1}(U_0)$
by $G$-invariant sections, hence $\sO_\P(g)$ is globally
generated by sections $s_1,\ldots,s_\ell\in H^0(\P,\sO_\P(g))$
over $\pi^{-1}(U_0)$.  By the Nullstellensatz, there exists some
$\beta'$ such that
$$
S^{\beta'}\Big(\bigoplus_{i=1}^\ell \sO_\P\cdot s_i\Big) \to
S^{g\cdot\beta'}(\sF)=\pi_* \sO_\P(g\cdot\beta')
$$
is surjective over $U_0$.
\end{proof}

The main result of this article says, that the existence of
smooth families $F:V\to U$ with ${\rm Var}(f) > 0$ is only
possible if $U$ carries multi-differential forms with logarithmic
singularities at infinity.

\begin{theorem}\label{0.3}
Let $Y$ be a projective manifold, S a reduced normal crossing
divisor, and let $f:V\to U=Y\setminus S$ be a smooth family
of $n$-dimensional projective varieties.
\begin{enumerate}
\item[i)] If $\omega_{V/U}$ is $f$-ample, then for some
$m >0$ the sheaf $S^m(\Omega^1_Y(\log S))$ contains an invertible sheaf
$\sA$ of Kodaira dimension $\kappa(\sA)\geq {\rm Var}(f)$.
\item[ii)] If $\omega_{V/U}$ is $f$-ample and
${\rm Var}(f)=\dim(Y)$, then for some $0< m\leq n$ the sheaf
$S^m(\Omega_Y^1(\log S))$ contains a big coherent subsheaf $\sP$.
\item[iii)] If $\omega_{V/U}$ is $f$-semi-ample and
${\rm Var}(f)=\dim(Y)$, then for some $m>0$ the sheaf
$S^m(\Omega_Y^1(\log S))$ contains a big coherent subsheaf $\sP$.
\item[iv)] Moreover under the assumptions made in iii)
there exists a non-singular finite covering $\psi:Y'\to Y$
and, for some $0< m \leq n$, a big coherent subsheaf $\sP'$ of
$\psi^*S^m(\Omega^1_Y(\log S))$.
\end{enumerate}
\end{theorem}

Before giving a guideline to the proof of \ref{0.3}, let us
discuss further properties of the sheaf of one forms on $U$, we
hope to be true.

\begin{problem} \label{0.2}
Let $Y$ be a projective manifold, $S$ a reduced normal crossing
divisor, and $U = Y \setminus S$. Let $\varphi : U \to M_h$ be a
morphism, induced by a family $f:V\to U$.
Assume that the family $f:V\to U$ induces an
\'etale map to the moduli stack, or in down to earth terms, that
the induced Kodaira Spencer map
$$
T_U \>>> R^1f_* T_{V/U}
$$
is injective and locally split.
\begin{enumerate}
\item[a)] Is $\Omega^{1}_{Y} (\log S)$ weakly positive,
or perhaps even weakly positive over $U$?
\item[b)] Is $\det (\Omega^{1}_{Y} (\log S ))=
\omega_Y(S)$ big?
\item[c)] Are there conditions on $\Omega^{1}_{F}$, for a
general fibre $F$ of $f$, which imply that $\Omega^{1}_{Y} (\log
 S)$ is big?
\end{enumerate}
\end{problem}

As we will see in \ref{4.1}, theorem \ref{0.3} implies that
the bigness in \ref{0.2}, b), follows from the weak positivity
in a).

There is hope, that the questions a) and b), which have been
raised by the first named author some time ago, will have an
affirmative answer. In particular they have been verified
by the second named author \cite{Zuo}, under the additional
assumption that the local Torelli theorem holds true for the
general fibre $F$ of $f$. The Brody hyperbolicity of moduli
stacks of canonically polarized manifolds, shown in \cite{V-Z2},
the results of Kov\'acs, and the content of this paper
strengthen this hope. As S. Kov\'acs told us, for certain divisors
$S$ in $Y=\P^N$, \ref{0.2}, a), holds true.

For moduli spaces of curves the sheaf $\Omega^{1}_{Y} (\log
S)$ is ample with respect to $U$. This implies
that morphisms $\pi:C_0\to U$ are rigid (see \ref{5.6}).
In the higher dimensional case the latter obviously does not
hold true (see \ref{5.3}), and problem c) asks for conditions implying
rigidity.

There is no evidence for the existence of a
reasonable condition in c). One could hope that
``$\Omega^{1}_{F}$ ample'' or ``$\Omega^{1}_{F}$ big'' will
work. At least, this excludes the obvious counter examples for
the ampleness of $\Omega^1_Y(\log S)$, discussed in \ref{5.3}.
For a non-isotrivial smooth family $V\to U$ of varieties with
$\Omega^1_F$ ample, the restriction of $\Omega^1_V$ to $F$
is big, an observation which for families of curves goes back to
H. Grauert \cite{Gra}. The problem \ref{0.1}, c), expresses our hope
that such properties of global multi-differential forms on the
general fibre could be mirrored in global properties of moduli spaces.\\

\noindent {\bf Notations.} To prove \ref{0.3} we start by choosing
any non-singular projective compactification $X$ of $V$, with
$\Delta= X\setminus V$ a normal crossing divisor, such that $V\to
U$ extends to a morphism $f:X\to Y$. For the proof of \ref{0.3} we
are allowed to replace $Y$ by any blowing up, if the pullback of
$S$ remains a normal crossing divisor. Moreover, as explained in
the beginning of the next section, we may replace $Y$ by the
complement of a codimension two subscheme, and $X$ by the
corresponding preimage, hence to work with partial
compactifications, as defined in \ref{1.1}. By abuse of notations,
such a partial compactification will again be denoted by $f:X\to
Y$.

In the course of the argument we will be forced to
replace the morphism $f$ by some fibred product. We will try
to keep the following notations. A morphism $f':X'\to Y'$ will
denote a pullback of $f$ under a morphism, usually dominant,
$Y' \to Y$, or a desingularization of such a pullback. The
smooth parts will be denoted by $V\to U$ and $V'\to U'$,
respectively. $f^r:X^r\to Y$ will denote the family obtained
as the $r$-fold fibred product over $Y$, and $f^{(r)}:X^{(r)}\to
Y$ will be obtained as a desingularization of $X^r$.
Usually $U_0$ will denote an open dense subscheme of $Y$,
and $\tilde{U}$ will be a blowing up of $U$.

At several places we need in addition some auxiliary
constructions. In section \ref{1} this will be a family
$g:Z \to Y'$, dominating birationally $X'\to Y'$ and a specific
model $g':Z'\to Y'$. For curves $C$ mapping to $Y$, the
desingularization of the induced family will be $h:W\to C$, where
again some  ${}'$ is added whenever we have to consider a
pullback family over some covering $C'$ of $C$.

Finally, in section \ref{3} $h:W \to Y$ will be a blowing up of
$X \to Y$ and $g:Z \to Y$ will be obtained as the
desingularization of a finite covering of $W$.\\
\ \\
{\bf Outline of the proof of \ref{0.3}.}
Let us start with \ref{0.3}, iii). In section \ref{2}
we will formulate and recall certain positivity properties of
direct image sheaves. In particular, by  \cite{Vie1}
the assumptions in iii) imply that $\det(f_*\omega_{X/Y}^\nu)$ is big,
for some $\nu \gg 2$. This in turn implies that the sheaf
$f_*\omega_{X/Y}^\nu$ is big. In \ref{2.10}
we will extend this result to the slightly smaller sheaf
$f_*(\Omega^n_{X/Y}(\log \Delta))^\nu.$ Hence for an ample
invertible sheaf $\sA$ on $Y$ and for some $\mu \gg 1$
the sheaf $S^\mu(f_*(\Omega^n_{X/Y}(\log \Delta))^\nu)\otimes \sA^{-1}$
will be globally generated over an open dense subset $U_0$.
Replacing $f:X\to Y$ by a partial compactification we will
assume this sheaf to be locally free. Then, for $\nu$
sufficiently large and divisible, $\Omega^n_{X/Y}(\log
\Delta))^{\nu\mu}\otimes f^*\sA^{-1}$ will be globally generated
over $f^{-1}(U_0)$, for some $U_0\neq \emptyset$.

This statement is not strong enough. We will need that
\begin{equation}\label{ass}
\Omega^n_{X/Y}(\log \Delta)^\nu \otimes f^* \sA^{-\nu}\mbox{ \  is globally
generated over \ } f^{-1}(U_0),
\end{equation}
for some $\nu \gg 2$ and for some ample invertible sheaf $\sA$.
To this aim we replace in \ref{2.10} the original morphism $f:X
\to Y$ by the $r$-th fibred product $f^{(r)}:X^{(r)} \to Y$,
for some $r\gg 1$.

The condition (\ref{ass}) will reappear in section \ref{3} in (\ref{3.3.1}).
There we study certain Higgs bundles $\bigoplus F^{p,q}$.
(\ref{3.3.1}) will allow to show, that $\sA\otimes \bigoplus
F^{p,q}$ is contained in a Higgs bundles induced by a variation
of Hodge structures. The latter is coming from a finite cyclic
covering $Z$ of $X$. The negativity theorem in \cite{Zuo} will
finish the proof of \ref{0.3} iii).

For each of the other cases in \ref{0.3} we need some additional
constructions, most of which are discussed in section \ref{1}. In
ii) and iv) (needed to prove \ref{cor2}, a) we
have to bound $m$ by the fibre dimension, hence we are not allowed
to replace $f:X\to Y$ by the fibre product
$f^{(r)}:X^{(r)} \to Y$. Instead we choose a
suitable covering $Y'\to Y$, in such a way that the assumption
(\ref{3.3.1}) holds true over the covering. For i) we have to
present the family $X'\to Y'$ as the pullback of a
family of maximal variation.

In section \ref{1} we also recall the weak semi-stable reduction
theorem due to Abramovich and Karu \cite{A-K} and some of its
consequences. In particular it will allow to construct a
generically finite dominant morphism $Y'\to Y$ such that for a
desingularization of the pullback family $f':X'\to Y'$ the sheaves
$\bigotimes^\mu f'_* \omega^\nu_{X'/Y'}$ are reflexive. This fact
was used in \cite{V-Z2} in the proof of \ref{2.10}. As mentioned,
we can restrict ourselves to partial compactifications $f:X\to Y$,
and repeating the arguments from \cite{V-Z2} in this case, we
would not really need the weakly semi-stable reduction. However,
the proof of the finiteness theorem \ref{5.2}
is based on this method.\\

As in \cite{V-Z} it should be sufficient in \ref{0.3}, iii) and
iv), to require that the fibres $F$ of $f$ are of general type,
or in the case $0 \leq \kappa(F) < \dim(F)$ that $F$ is
birational to some $F'$ with $\omega_{F'}$ semi-ample.
We do not include this, since the existence of relative
base loci make the notations even more confusing than they are
in the present version. However, comparing the arguments in
\cite{V-Z}, \S 3, with the ones used here, it should not be too
difficult to work out the details.

\section{Mild morphisms}\label{1}

As explained at the end of the last section it will be
convenient, although not really necessary, to use for the
proof of \ref{0.3} some of the results and constructions
contained in \cite{V-Z2}, in particular the weak semi-stable
reduction theorem due to Abramovich and Karu \cite{A-K}. It will
allow us to formulate the strong positivity theorem \ref{2.10}
for product families, shown in \cite{V-Z2}, 4.1., and it will be
used in the proof of the boundedness of the functor of
homomorphism in \ref{5.2}. We will use it again to reduce the
proof of \ref{0.3}, i), to the case of maximal
variation, although this part could easily be done without the weak
semi-stable reduction. We also recall Kawamata's covering
construction. The content of this section will be needed in
the proof of parts i), ii) and iv) of \ref{0.3}, but not for
iii).

\begin{definition}\label{1.1}\ \ \
\begin{enumerate}
\item[a)] Given a family $\tilde{V}\to \tilde{U}$ we will call
$V\to U$ a birational model if there exist compatible birational
morphisms $\tau:U\to \tilde{U}$ and
$\tau':V\to \tilde{V}\times_{\tilde{U}}U$. If we underline that
$U$ and $\tilde{U}$ coincide, we want $\tau$ to be the identity.
If $\tilde{V}\to \tilde{U}$ is smooth, we call $V\to U$ a smooth
birational model, if $\tau'$ is an isomorphism.
\item[b)] If $V\to U$ is a smooth projective family of
quasi-projective manifolds, we call $f:X\to Y$ a partial
compactification, if
\begin{enumerate}
\item[i)] $X$ and $Y$ are quasi-projective manifolds, and $U\subset Y$.
\item[ii)] $Y$ has a non-singular projective compactification
$\bar{Y}$ such that $S$ extends to a normal crossing
divisor and such that ${\rm codim}(\bar{Y}\setminus Y) \geq 2$.
\item[iii)] $f$ is a projective morphism and $f^{-1}(U)\to U$
coincides with $V\to U$.
\item[iv)] $S=Y\setminus U$, and $\Delta=f^*S$ are normal
crossing divisors.
\end{enumerate}
\item[c)] We say that a partial compactification $f:X\to Y$ is a
good partial compactification if the condition iv) in b) is
replaced by
\begin{enumerate}
\item[iv)] $f$ is flat, $S=Y\setminus U$ is a smooth divisor, and
$\Delta=f^*S$ is a relative normal crossing divisor, i.e. a
normal crossing divisor whose components, and all their intersections
are smooth over components of $S$.
\end{enumerate}
\item[d)] The good partial compactification $f:X\to Y$ is
semi-stable, if in c), iv), the divisor $f^*S$ is reduced.
\item[e)] An arbitrary partial compactification of $V\to U$ is
called semi-stable in codimension one, if it contains a
semi-stable good partial compactification.
\end{enumerate}
\end{definition}
\begin{remark}\label{kappa}
The second condition in b) or c) allows to talk about
invertible sheaves $\sA$ of positive Kodaira dimension on $Y$.
In fact, $\sA$ extends in a unique way to an invertible sheaf
$\bar{\sA}$ on $\bar{Y}$ and
$$
H^0(Y,\sA^\nu) = H^0(\bar{Y},{\bar{\sA}}^\nu),
$$
for all $\nu$. So we can write $\kappa(\sA):=\kappa(\bar{\sA})$,
in case $Y$ allows a compactification satisfying ii).
If $\tau: \bar{Y}' \to \bar{Y}$ is a blowing up with centers
in $\bar{Y}\setminus Y$, and if $\bar{\sA}'$ is an extension
of $\sA$ to $\bar{Y}'$, then $\kappa(\sA) \geq
\kappa(\bar{\sA}')$.

In a similar way, one finds a coherent sheaf $\sF$ on $Y$ to be semi-ample
with respect to $U_0\subset Y$ (or weakly positive over $U_0$), if
and only if its extension to $\bar{Y}$ has the same property.
\end{remark}

Kawamata's covering construction will be used frequently
throughout this article. First of all, it allows the semi-stable
reduction in codimension one, and secondly it allows to take
roots out of effective divisors.
\begin{lemma}\label{1.2} \ \ \
\begin{enumerate}
\item[a)] Let $Y$ be a quasi-projective manifold, $S$ a normal
crossing divisor, and let $\sA$ be an invertible sheaf, globally
generated over $Y$. Then for all $\mu$ there exists a
non-singular finite covering $\psi:Y'\to Y$ whose discriminant
$\Delta(Y'/Y)$ does not contain components of $S$, such that
$\psi^*(S+\Delta(Y'/Y))$ is a normal crossing divisor, and such that
$\psi^*\sA=\sO_{Y'}(\mu\cdot A')$ for some reduced
non-singular divisor $A'$ on $Y'$.
\item[b)] Let $f:X\to Y$ be a partial compactification of a
smooth family $V\to U$. Then there exists a non-singular finite
covering $\psi:Y'\to Y$, and a desingularization $\psi':X'\to
X\times_YY'$ such that the induced family $f':X'\to Y'$ is
semi-stable in codimension one.
\end{enumerate}
\end{lemma}
\begin{proof}
Given positive integers $\epsilon_i$ for all components
$S_i$ of $S$, Kawamata constructed a finite non-singular
covering $\psi:Y'\to Y$ (see \cite{Vie}, 2.3), with
$\psi^*(S+\Delta(Y'/Y))$ a normal crossing divisor,
such that all components of $\psi^*S_i$ are ramified of order
exactly $\epsilon_i$.

In a) we choose $A$ to be the zero-divisor of a general section
of $\sA$, and we apply Kawamata's construction to $S+A$, where
the $\epsilon_i$ are one for the components of $S$ and where the
prescribed ramification index for $A$ is $\mu$.

In b) the semi-stable reduction theorem for families over
curves (see \cite{KKMS}) allows to choose the $\epsilon_i$ such
that the family $f':X'\to Y'$ is semi-stable in codimension one.
\end{proof}

Unfortunately in \ref{1.2}, b), one has little control on
the structure of $f'$ over the singularities of $S$. Here
the weak semi-stable reduction theorem will be of help.
The pullback of a weakly semi-stable morphism under a dominant
morphism of manifolds is no longer weakly semi-stable. However
some of the properties of a weakly semi-stable morphism survive. Those
are collected in the following definition, due again to Abramovich and
Karu \cite{A-K}.

\begin{definition} \label{1.3}
A projective morphism $g' : Z' \to Y'$ between quasi-projective
varieties is called mild, if
\begin{enumerate}
\item[a)] $g'$ is flat, Gorenstein with reduced fibres.
\item[b)] $Y'$ is non-singular and $Z'$ normal with at most rational
singularities.
\item[c)] Given a dominant morphism $Y'_1 \to Y'$ where $Y'_1$ has at
most rational Gorenstein singularities, $Z' \times_{Y'} Y'_1$ is
normal with at most rational singularities.
\item[d)] Let $Y'_0$ be an open subvariety of $Y'$, with
${g'}^{-1}(Y'_0)\to Y'_0$ smooth.
Given a non-singular curve $C'$ and a morphism
$\pi:C' \to Y'$ whose image meets $Y'_0$, the fibred product
$Z'\times_{Y'} C'$ is normal, Gorenstein with at most rational
singularities.
\end{enumerate}
\end{definition}

The mildness of $g'$ is, more or less by definition, compatible with
pullback. Let us rephrase three of the properties shown in \cite{V-Z2}, 2.2.

\begin{lemma} \label{1.4} Let $Z$ and $Y$ be quasi-projective
manifolds, $g : Z \to Y'$ be a projective, birational
to a projective mild morphism $g':Z'\to Y'$. Then one has:
\begin{enumerate}
\item[i)] For all $\nu \geq 1$ the sheaf $g_* \omega^{\nu}_{Z/Y'}$ is
reflexive and isomorphic to $g'_*\omega_{Z'/Y'}^\nu$.
\item[ii)] If $\gamma:Y'' \to Y'$ is a dominant morphism between
quasi-projective manifolds, then the morphism
$pr_2:Z \times_{Y'} Y'' \to Y''$ is birational to a
projective mild morphism to $Y''$.
\item[iii)] Let $Z^{(r)}$ be a desingularization of the $r$-fold
fibre product $Z\times_{Y'} \cdots \times_{Y'}Z$. Then the
induced morphism $Z^{(r)}\to Y'$ is birational to a projective
mild morphism over $Y'$.
\end{enumerate}
\end{lemma}

One consequence of the weakly semi-stable reduction says
that, changing the birational model of a morphism, one always
finds a finite cover of the base such that the pullback is
birational to a mild morphism (see \cite{V-Z2}, 2.3).
\begin{lemma}\label{1.5}
Let $V\to U$ be a smooth family of projective
varieties. Then there exists a quasi-projective manifold
$\tilde{U}$ and a smooth birational model $\tilde{V}\to
\tilde{U}$, non-singular projective compactifications $Y$ of
$\tilde{U}$ and $X$ of $\tilde{V}$, with
$S=Y\setminus \tilde{U}$ and $\Delta=X\setminus \tilde{V}$ normal crossing
divisors, and a diagram of projective morphisms
\begin{equation}\label{1.5.1}
\begin{CDS}
X \< \psi' << X' \< \sigma << Z \\
\V V f V \novarr \V V f' V \SW W g W \\
Y \< \psi << Y'\novarr
\end{CDS}
\end{equation}
with:
\begin{enumerate}
\item[a)] $Y'$ and $Z$ are non-singular, $X'$ is the
normalization of $X\times_YY'$, and $\sigma$ is a desingularization.
\item[b)] $g$ is birational to a mild morphism $g':Z'\to Y'$.
\item[c)] For all $\nu >0$ the sheaf
$g_*\omega^{\nu}_{Z/Y'}$ is reflexive and there exists an injection
$$
g_*\omega_{Z/Y'}^\nu \>>> \psi^*f_*\omega_{X/Y}^\nu.
$$
\item[d)] For some positive integer $N_\nu$, and for some
invertible sheaf $\lambda_\nu$ on $Y$
$$
\det(g_*\omega_{Z/Y'}^\nu)^{N_\nu}=\psi^*\lambda_\nu.
$$
\item[e)] Moreover, if $\bar{Y}$ is a given projective compactification of $U$,
we can assume that there is a birational morphism $Y\to \bar{Y}.$ 
\end{enumerate}
\end{lemma}
This diagram has been constructed in \cite{V-Z2}, \S 2. Let us just
recall that the reflexivity of $g_*\omega^{\nu}_{Z/Y'}$ in c) is
a consequence of b), using \ref{1.4}, i).\qed \\

Unfortunately the way it is constructed, the mild model
$g':Z'\to Y'$ might not be smooth over $\psi^{-1}(\tilde{U})$.
Moreover even in case $U$ is non-singular one has to allow
blowing ups $\tau:\tilde{U}\to U$. Hence starting from $V\to U$ we can
only say that $U$ contains some ``good'' open dense subset $U_g$
for which $\tau$ is an isomorphism between
$\tilde{U}_g:=\tau^{-1}(U_g)$ and $U_g$, and for which
$$
{g'}^{-1}\psi^{-1}(\tilde{U}_g) \>>> \psi^{-1}(\tilde{U}_g)
$$
is smooth. The next construction will be needed in the proof of
\ref{5.2}.
\begin{corollary}\label{1.6}
Let $C$ be a non-singular projective curve, $C_0\subset C$ an
open dense subset, and let $\pi_0 : C_0 \to U$
be a morphism with $\pi_0(C_0)\cap U_g\neq \emptyset$. Hence
there is a lifting of $\pi_0$ to $\tilde{U}$ and an extension $\pi:C\to
Y$ with $\pi_0 = \tau\circ\pi|_{C_0}$. Let $h: W\to C$ be the
family obtained by desingularizing the main component of the
normalization of $X\times_YC$, and let $\lambda_\nu$ and $N_\nu$
be as in \ref{1.5}, d). Then
$$\deg(\pi^*\lambda_\nu) \leq N_\nu\cdot
\deg(\det(h_*\omega_{W/C}^\nu)). $$
\end{corollary}
\begin{proof}
Let $\rho:C' \to C$ be a finite morphism of
non-singular curves such that $\pi$ lifts to $\pi':C'\to Y'$,
and let $h': W'\to C'$ be the family obtained by
desingularizing the main component of the normalization of
$X'\times_{Y'}C'$. By condition d) in the definition of a mild
morphism, and by the choice of $U_g$, the family $h'$ has
$Z'\times_{Y'}C'$ as a mild model. Applying \ref{1.5}, c),
to $h$ and $h'$, and using \ref{1.5}, d), we find
\begin{gather*}
\deg(\rho)\cdot\deg(\pi^*\lambda_\nu)=
{N_\nu}\cdot\deg({\pi'}^*g_*\omega_{Z/Y'}^\nu) \mbox{ \ \ \ and}\\
\deg(h'_*\omega_{W'/C'}^\nu) \leq \deg(\rho)\cdot\deg(h_*\omega_{W/C}^\nu).
\end{gather*}
Moreover, by \ref{1.5}, c), and by base change, one
obtains a morphism of sheaves
\begin{equation}\label{1.6.1}
{\pi'}^*g_*\omega_{Z/Y'}^\nu \simeq
{\pi'}^*g'_*\omega_{Z'/Y'}^\nu \>>>
pr_{2*}\omega_{Z'\times_{Y'}C'/C'}^\nu \simeq
h'_*\omega_{W'/C'}^\nu,
\end{equation}
which is an isomorphism over some open dense subset.
Let $r$ denote the rank of those sheaves.

It remains to show, that (\ref{1.6.1}) induces an injection
from ${\pi'}^*\det(g'_*\omega_{Z'/Y'}^\nu)$ into
$\det(pr_{2*}\omega_{Z'\times_{Y'}C'/C'}^\nu)$.
To this aim, as in \ref{1.4}, iii), let ``${}^{(r)}$'' stand for
``taking a desingularization of the $r$-th fibre product''.
Then $g^{(r)}:Z^{(r)}\to Y'$ is again birational to a mild
morphism over $Y'$. As in \cite{V-Z2}, 4.1.1, flat base change
and the projection formula give isomorphisms
$$
{h'}^{(r)}_*\omega_{{W'}^{(r)}/C'}^\nu \simeq \bigotimes^r
h'_*\omega_{W'/C'}^\nu
\mbox{ \ \ \ and \ \ \ }
g^{(r)}_*\omega_{Z^{(r)}/Y'}^\nu \simeq \bigotimes^r
g_*\omega_{Z/Y'}^\nu,
$$
the second one outside of a codimension two subscheme.
Since both sheaves are reflexive, the latter extends to $Y'$.

The injection (\ref{1.6.1}),
applied to $g^{(r)}$ and $h^{(r)}$ induces an injective morphism
\begin{equation}\label{1.6.2}
{\pi'}^*\bigotimes^r g_*\omega_{Z/Y'}^\nu
\>>> \bigotimes^r h'_*\omega_{W'/C'}^\nu.
\end{equation}
The left hand side contains
${\pi'}^*\det(g'_*\omega_{Z'/Y'}^\nu)$ as direct factor, whereas
the righthand contains
$\det(pr_{2*}\omega_{Z'\times_{Y'}C'/C'}^\nu)$, and we obtain the
injection for the determinant sheaves as well.
\end{proof}
The construction of (\ref{1.5.1}) in \cite{V-Z2} will be used
to construct a second diagram (\ref{1.7.1}). There we do not insist on the
projectivity of the base spaces, and we allow ourselves to
work with good partial compactifications of an open subfamily of
the given one. This quite technical construction will be needed
to proof \ref{0.3}, i).

\begin{lemma}\label{1.7}
Let $V\to U$ be a smooth family of canonically polarized
manifolds. Let $\bar{Y}$ and $\bar{X}$ be non-singular projective
compactifications of $U$ and $V$ such that both, $\bar{Y}\setminus
U$ and $\bar{X}\setminus V$, are normal crossing divisors and
such that $V\to U$ extends to $\bar{f}:\bar{X}\to \bar{Y}$.
Blowing up $\bar{Y}$ and $\bar{X}$, if necessary, one finds
an open subscheme $Y$ in $\bar{Y}$ with ${\rm
codim}(\bar{Y}\setminus Y) \geq 2$ such that the
restriction $f:Y\to X$ is a good partial compactification of a
smooth birational model of $V\to U$, and one finds a diagram of
morphisms between quasi-projective manifolds
\begin{equation}\label{1.7.1}
\begin{CDS}
X \< \psi' << X' \< \sigma << Z \> \eta' >> Z^{\#}\\
\V V f V \novarr \V V f' V \SW W g W \novarr \SW W g^{\#} W \\
Y \< \psi << Y' \> \eta >> Y^{\#}
\end{CDS}
\end{equation}
with:
\begin{enumerate}
\item[a)] $g^\#$ is a projective morphism, birational to a mild
projective morphism ${g^\#}':{Z^\#}' \to Y^\#$.
\item[b)] $g^{\#}$ is semi-stable in codimension one.
\item[c)] $Y^\#$ is projective, $\eta$ is dominant and smooth,
$\eta'$ factors through a birational morphism $Z\to Z^\#\times_{Y^\#}Y'$,
and $\psi$ is finite.
\item[d)] $X'$ is the normalization of $X\times_YY'$ and
$\sigma$ is a blowing up with center in
${f'}^{-1}\psi^{-1}(S')$. In particular $f'$ and $g'$ are
projective.
\item[e)] Let $U^{\#}$ be the largest subscheme of $Y^\#$ with
$$
V^{\#}={g^\#}^{-1}(U^\#) \>>> U^{\#}
$$
smooth. Then $\psi^{-1}(U) \subset \eta^{-1}(U^\#)$,
and $U^\#$ is generically finite over $M_h$.
\item[f)] For all $\nu >0$ there are isomorphisms
$$
g_*\omega^{\nu}_{Z/Y'} \simeq \eta^* g^{\#}_* \omega^\nu_{Z^{\#}/Y^{\#}}
\mbox{ \ \ \ and \ \ \ } \det(g_*\omega_{Z/Y'}^\nu) \simeq
\eta^*\det(g^{\#}_*\omega_{Z^{\#}/Y^{\#}}^\nu).
$$
\item[g)] For all $\nu > 0$ there exists an injection
$$
g_*(\omega_{Z/Y'}^\nu) \>>> \psi^*f_*(\omega_{X/Y}^\nu).
$$
For some positive integer $N_\nu$, and for some
invertible sheaf $\lambda_\nu$ on $Y$
$$
\det(g_*(\omega_{Z/Y'}^\nu))^{N_\nu}=\psi^*\lambda_\nu.
$$
\end{enumerate}
\end{lemma}
\begin{proof}
It remains to verify, that the construction given in
\cite{V-Z2}, \S 2, for (\ref{1.5.1}) can be modified
to guaranty the condition e) along with the others.

Let $\varphi:U\to M_h$ be the induced morphism to the moduli scheme.
Seshadri and Koll\'ar constructed a finite Galois cover of the
moduli space which is induced by a family (see \cite{Vie}, 9.25,
for example). Hence there exists some manifold ${U}^{\#}$,
generically finite over the closure of $\varphi(U)$ such that
the morphism ${U}^{\#} \to  M_h$ is induced by a family $V^{\#} \to U^{\#}$.
By \cite{V-Z2}, 2.3, blowing up $U^\#$, if necessary, we find a projective
compactification $Y^\#$ of $U^\#$ and a covering ${Y^\#}'$, such
that
$$
V^\#\times_{{Y^\#}}{Y^\#}' \>>> {Y^\#}'
$$
is birational to a projective mild morphism over ${Y^\#}'$.
Replacing $U^\#$ by some generically finite cover, we can assume
that $V^\# \to U^\#$ has such a model already over $Y^\#$.

Next let $Y'$ be any variety, generically finite over $\bar{Y}$,
for which there exists a morphism $\eta:Y'\to Y^\#$.
By \ref{1.4}, ii), we are allowed to replace $Y^\#$ by any
manifold, generically finite over $Y^\#$, without loosing the
mild birational model. Doing so, we can assume the fibres of
$Y'\to Y^\#$ to be connected. Replacing
$Y'$ by some blowing up, we may assume that for some
non-singular blowing up $Y\to \bar{Y}$ the morphism $Y'\to
\bar{Y}$ factors through a finite morphism $Y'\to Y$.

Next choose a blowing up ${Y^\#}'
\to Y^\#$ such that the main component of
$Y'\times_{Y^\#}{Y^\#}'$ is flat over $Y^\#$, and $Y''$ to be a
desingularization. Hence changing notations again,
and dropping one prime, we can assume that
the image of the largest reduced divisor $E$ in $Y'$ with ${\rm
codim}(\eta(E)) \geq 2$ maps to a subscheme of $Y$ of
codimension larger that or equal to $2$. This remains true, if
we replace $Y^\#$ and $Y'$ by finite coverings.
Applying \ref{1.2}, b), to $Y' \to Y^\#$, provides us with
non-singular covering of $Y^\#$ such that a desingularization of
the pullback of $Y'\to Y^\#$ is semi-stable in codimension one.
Again, this remains true if we replace $Y^\#$ by a larger
covering, and using \ref{1.2}, b), a second time, now for
$Z^\# \to Y^\#$, we can assume that this morphism is as well
semi-stable in codimension one.

Up to now, we succeeded to find the manifolds in (\ref{1.7.1})
such that a) and b) hold true. In c), the projectivity of $Y^\#$
and the dominance of $Y'$ over $Y^\#$ follow from the
construction. For the divisor $E$ in $Y'$ considered above,
we replace $Y$ by $Y\setminus\psi(E)$ and $Y'$ by $Y'\setminus
E$, and of course $X$, $X'$ and $Z$ by the corresponding preimages.
Then the non-equidimensional locus of $\eta$ in $Y'$ will be of
codimension larger than or equal to two. $\psi$ is
generically finite, by construction, hence finite over the
complement of a codimension two subscheme of $Y$.
Replacing again $Y$ by the complement of codimension two subscheme, we
can assume $\eta$ to be equidimensional, hence flat, and
$\psi$ to be finite. The morphism $\eta$ has reduced fibres over
general points of divisors in $Y^\#$, hence it is smooth outside
a codimension two subset of $Y'$, and replacing $Y$ by the
complement of its image, we achieved c).

Since $V\to U$ is smooth, the pullback of $X\to Y$ to $Y'$ is
smooth outside of $\psi^{-1}(S)$. Moreover the induced morphism
to the moduli scheme $M_h$ factors through an open subset of
$Y^\#$. Since by construction $U^\#$ is proper over its image in
$M_h$, the image of $\psi^{-1}(U)$ lies in $U^\#$ and we obtain
d) and e).

For f) remark that the pullback $Z'\to Y'$ of the mild
projective morphism $g^\#:{Z^\#}' \to Y^\#$ to $Y'$ is again
mild, and birational to $Z\to Y'$. By flat base change,
$$
g'_*\omega^\nu_{Z'/Y'} \simeq \eta^*{g^\#}'_*\omega_{{Z^\#}'/Y^\#}^\nu.
$$
Since $Z'$ and ${Z^\#}'$ are normal with at most rational
Gorenstein singularities, we obtain (as in \cite{V-Z2}, 2.3)
that the sheaf on the right hand side is
$g_*\omega^\nu_{Z/Y'}$ whereas the one on the right hand side is
$\eta^*{g^\#}_*\omega_{{Z^\#}/Y^\#}^\nu$.

g) coincides with \ref{1.5}, d), and it has been verified in
\cite{V-Z2}, 2.4. as a consequence of the existence of a mild
model for $g$ over $Y'$.
\end{proof}

\section{Positivity and ampleness}\label{2}

Next we will recall positivity theorems, due to Fujita, Kawamata,
Koll\'ar and the first named author. Most of the content of this
section is well known, or easily follows from known results.

As in \ref{0.3} we will assume throughout
this section that $U$ is the complement of a normal crossing
divisor $\bar{S}$ in a manifold $\bar{Y}$, and that there is a
smooth family $V\to U$ with $\omega_{V/U}$ relative semi-ample.
Leaving out a codimension two subset in $\bar{Y}$ we
find a good partial compactification $f:X\to Y$, as defined in
\ref{1.1}.

For an effective $\Q$-divisor $D\in {\rm Div}(X)$ the integral
part $[D]$ is the largest divisor with $[D] \leq D$.
For an effective divisor $\Gamma $ on $X$, and for
$N\in \N-\{0\}$ the algebraic multiplier sheaf is
$$
\omega_{X/Y} \Big\{ \frac{-\Gamma}{N} \Big\} = \psi_* \Big(
\omega_{T/Y} \Big( - \Big[ \frac{\Gamma'}{N} \Big] \Big)
\Big)
$$
where $\psi : T \to X$ is any blowing up with $\Gamma' = \psi^*
\Gamma$ a normal crossing divisor (see for example \cite{E-V},
7.4, or \cite{Vie}, section 5.3).

Let $F$ be a non-singular fibre of $f$.
Using the definition given above for $F$, instead of $X$, and for a divisor
$\Pi$ on $F$, one defines
$$
e (\Pi) = {\rm Min} \Big\{ N \in \N \setminus  \{ 0 \} ; \ \omega_F
\Big\{ \frac{- \Pi}{N} \Big\} = \omega_F \Big\} .
$$
By \cite{E-V} or \cite{Vie}, section 5.4, $e(\Gamma |_{F})$ is
upper semi-continuous, and there exists a neighborhood $V_0$ of $F$
with $e (\Gamma |_{V_0} ) \leq e (\Gamma |_{F})$. If $\sL$ is an
invertible sheaf on $F$, with $H^0(F,\sL)\neq 0$, one defines
$$
e (\sL) = {\rm Max} \big\{e(\Pi) ; \ \Pi \mbox{ an effective
divisor and } \sO_F(\Pi)=\sL \big\} .
$$
\begin{proposition}[\cite{V-Z2}, 3.3] \label{2.3}
Let $\sL$ be an invertible sheaf, let $\Gamma$ be a divisor on
$X$, and let $\sF$ be a coherent sheaf on $Y$. Assume that, for
some $N >0$ and for some open dense subscheme $U_0$ of $U$, the
following conditions hold true:
\begin{enumerate}
\item[a)] $\sF$ is weakly positive over $U_0$ (in particular $\sF
|_{U_0}$ is locally free).
\item[b)] There exists a morphism $f^* \sF \to \sL^N  ( -
\Gamma)$, surjective over $f^{-1}(U_0)$.
\item[c)] None of the fibres $F$ of $f: V_0=f^{-1}(U_0) \to U_0$
is contained in $\Gamma$, and for all of them
$$
e (\Gamma |_{F}) \leq N.
$$
\end{enumerate}
Then $f_* (\sL \otimes \omega_{X/Y})$ is weakly positive over
$U_0$.
\end{proposition}
As mentioned in \cite{V-Z2}, 3.8, the arguments used in
\cite{Vie}, 2.45, carry over to give a simple proof of the
following, as a corollary of \ref{2.3}.
\begin{corollary}[\cite{Vie3}, 3.7] \label{2.4}
$f_* \omega^{\nu}_{X/Y}$ is weakly positive over $U$.
\end{corollary}
In \cite{Vie2}, for families of canonically polarized
manifolds and in \cite{Kaw}, in general, one finds the
strong positivity theorem saying:
\begin{theorem}\label{2.5}
If $\omega_{V/U}$ is $f$-semi-ample, then for some $\eta$ sufficiently large
and divisible,
$$
\kappa(\det(f_*\omega_{X/Y}^\eta))\geq {\rm Var}(f).
$$
\end{theorem}

In case ${\rm Var}(f)=\dim(Y)$, \ref{2.5} implies that
$\det(f_*\omega_{X/Y}^\eta)$ is ample with respect to some open
dense subset $U_0$ of $Y$. If the general fibre of
$f$ is canonically polarized, and if the induced map
$\varphi:U\to M_h$ quasi-finite over its image, one can choose
$U_0=U$, as follows from the last part of the next proposition.

\begin{proposition} \label{2.6}
Assume that ${\rm Var}(f)=\dim(Y)$, and that $\omega_{V/U}$ is
$f$-semi-ample. Then:
\begin{enumerate}
\item[i)] The sheaf
$f_*\omega_{X/Y}^\nu$ is ample with respect to some open dense
subset $U_0$ of $Y$ for all $\nu > 1$ with
$f_*\omega_{X/Y}^\nu\neq 0$.
\item[ii)] If $B$ is an effective divisor, supported in $S$ then
for all $\nu$ sufficiently large and divisible, the sheaf $\sO_Y(-B)\otimes
f_*\omega_{X/Y}^\nu$ is ample with respect to some open dense
subset $U_0$.
\item[iii)] If the smooth fibres of $f$ are canonically
polarized and if the induced morphism $\varphi:U\to M_h$ is
quasi-finite over its image, then one can chose $U_0=U$ in i)
and ii).
\end{enumerate}
\end{proposition}
\begin{proof}
For iii) one uses a variant of \ref{2.5}, which has been shown
\cite{Vie3}, 1.19. It also follows from the obvious extension of
the ampleness criterion in \cite{Vie}, 4.33, to the case
``ample with respect to $U$'':
\begin{claim}\label{2.7}
Under the assumption made in \ref{2.6}, iii), for all $\eta$
sufficiently large and divisible, there exist positive integers
$a$, $b$ and $\mu$ such that
$$
\det (f_* \omega^{\mu\eta}_{X/Y})^a\otimes \det (f_*
\omega^{\eta}_{X/Y})^b
$$
 ample with respect to $U$. \qed
\end{claim}
Since we do not want to distinguish between the two cases i) and
iii) in \ref{0.3}, we choose $U_0=U$ in iii), and we allow
$a=\mu=0$ in case i). By \ref{2.7} and \ref{2.5}, respectively,
in both cases the sheaf $\det (f_*
\omega^{\mu\eta}_{X/Y})^a\otimes \det (f_*
\omega^{\eta}_{X/Y})^b$ is ample with respect to $U_0$.

By \cite{E-V}, \S 7, or \cite{Vie}, Section 5.4,
the number $e(\omega^{\mu\eta}_F)$ is bounded by some constant
$e$, for all smooth fibres of $f$. We will choose $e$ to be
divisible by $\eta$ and larger than $\mu\eta$.

Replacing $a$ and $b$ by some multiple, we may
assume that there exists a very ample sheaf $\sA$ and a morphism
$$
\sA \>>> \det (f_* \omega^{\mu\eta}_{X/Y})^a\otimes \det (f_*
\omega^{\eta}_{X/Y})^b
$$
which is an isomorphism over $U_0$, and that $b$ is divisible by
$\mu$.

By \ref{1.2}, a), there exists a non-singular covering $\psi:Y'
\to Y$ and an effective divisor $H$ with $\psi^*\sA
=\sO_{Y'}(e\cdot (\nu-1)\cdot H)$, and such that the
discriminant locus $\Delta(Y'/Y)$ does not contain any of the
components of $S$. Replacing $Y$ by a slightly smaller scheme,
we can assume that $\Delta(Y'/Y)\cap S = \emptyset$, hence $X'=
X\times_YY'$ is non-singular and by flat base change
$$
pr_{2*}\omega_{X'/Y'}^\sigma = \psi^*f_*\omega_{X/Y}^\sigma
$$
for all $\sigma$. The assumptions in \ref{2.6}, i), ii) or iii) remain
true for $pr_2:X'\to Y'$, and by \cite{Vie}, 2.16, it is
sufficient to show that the conclusions in \ref{2.6} hold true on $Y'$ for
$\psi^{-1}(U_0)$.

Dropping the primes, we will assume in the sequel that $\sA$ has
a section whose zero-divisor is $e\cdot (\nu-1)\cdot H$ for a non-singular
divisor $H$.

Let $r(\sigma)$ denote the rank of
$f_*\omega_{X/Y}^\sigma$. We choose
$$
r=r(\eta)\cdot \frac{b}{\mu}+r(\mu\eta)\cdot a,
$$
consider the $r$-fold fibre product
$$
f^r : X^r = X\times_Y X \ldots \times_Y X \>>> Y,
$$
and a desingularization $\delta : X^{(r)} \to X^r.$
Using flat base change, and the natural maps
$$
\sO_{X^r} \>>> \delta_*\sO_{X^{(r)}}\mbox{ \ \  and \ \ }
\delta_*\omega_{X^{(r)}} \>>> \omega_{X^{r}},
$$
one finds morphisms
\begin{gather}
\bigotimes^{r} f_* \omega^{\mu\eta}_{X/Y}
\>>> f^{(r)}_{*} \delta^*\omega^{\mu\eta}_{X^r/Y} \label{2.7.2} \ \
\mbox{ \ \ \ and}\\[.1cm]
f^{(r)}_{*} \delta^*(\omega_{X^r/Y}^{\nu-1}\otimes \omega_{X^{(r)} /Y}) \>>>
f^{r}_{*}\omega^{\nu}_{X^r/Y} =
\bigotimes^r f_* \omega^{\nu}_{X/Y}, \label{2.7.3}
\end{gather}
and both are isomorphism over $U$. We have natural maps
\begin{gather}\label{2.7.4}
\det(f_*\omega_{X/Y}^{\mu\eta}) \> \subset >> \bigotimes^{r(\mu\eta)}f_*
\omega_{X/Y}^{\mu\eta}  \mbox{ \ \ and}\\ \label{2.7.5}
\det(f_*\omega_{X/Y}^{\eta})^\mu \> \subset >> \bigotimes^{r(\eta)\cdot \mu}f_*
\omega_{X/Y}^{\eta} \>>> \bigotimes^{r(\eta)}f_*
\omega_{X/Y}^{\mu\eta},
\end{gather}
where the last morphism is the multiplication map. Hence we
obtain
\begin{multline*}
\sA=\sO_Y(e\cdot(\nu-1)\cdot H) \>>> \det (f_*
\omega^{\mu\eta}_{X/Y})^a\otimes \det (f_*
\omega^{\eta}_{X/Y})^b \>>>\\
\bigotimes^rf_*\omega_{X/Y}^{\mu\eta}\>>> f^{(r)}_{*}
\delta^*\omega^{\mu\eta}_{X^r/Y} .
\end{multline*}
Thereby the sheaf ${f^{(r)}}^* \sA$ is a subsheaf of
$\delta^*\omega_{X^r/Y}^{\mu\eta}$. Let $\Gamma$ be the zero divisor of
the corresponding section of
$$
f^{(r)*}\sA^{-1}\otimes\delta^*\omega^{\mu\eta}_{X^r/Y}
\mbox{ \ \ hence \ \ }\sO_{X^{(r)}}(-\Gamma)=
f^{(r)*}\sA\otimes\delta^*\omega^{-\mu\eta}_{X^r/Y} .
$$
For the sheaf
$$
\sM=\delta^* (\omega_{X^r/Y}\otimes
\sO_{X^r}(-f^{r*}H))
$$
one finds
$$
\sM^{e\cdot (\nu-1)}(-\Gamma) = \delta^*\omega_{X^r/Y}^{e\cdot
(\nu-1)}\otimes f^{(r)*}\sA^{-1}\otimes\sO_{X^{(r)}}(-\Gamma)
= \delta^*\omega_{X^r/Y}^{e\cdot (\nu-1) -\mu\eta}.
$$
By the assumption \ref{2.6}, i), and by the choice of $e$
we have a morphism
$$
f^*f_*\omega_{X/Y}^{e\cdot(\nu-1)-\mu\eta} \>>>
\omega_{X/Y}^{e\cdot(\nu-1)-\mu\eta},
$$
surjective over $f^{-1}(U_0)$. The sheaf
$$
\sF=\bigotimes^rf_*\omega_{X/Y}^{e\cdot(\nu-1)-\mu\eta}
$$
is weakly positive over $U_0$ and there is a morphism
$$
f^{(r)*}\sF \>>> \sM^{e\cdot (\nu-1)}(-\Gamma)
$$
surjective over $f^{-1}(U_0)$. Since the morphism of sheaves in
(\ref{2.7.4}), as well as the first one in (\ref{2.7.5}), split
locally over $U_0$ the divisor $\Gamma$ can not contain a fibre
$F$ of
$$
{f^{(r)}}^{-1}(U_0) \>>> U_0,
$$
and by \cite{E-V}, \S 7, or \cite{Vie}, 5.21, for those fibres
$$
e(\Gamma|_{F^r})\leq
e(\omega_{F^r}^{\mu\eta})=e(\omega_F^{\mu\eta})\leq e.
$$
Applying \ref{2.3} to $\sL=\sM^{\nu-1}$ one obtains the weak
positivity of the sheaf
$$
f^{(r)}_*(\sM^{\nu-1}\otimes \omega_{X^{(r)}/Y})=
f^{(r)}_*(\delta^*(\omega_{X^r/Y}^{\nu-1}\otimes
\omega_{X^{(r)}/Y}))\otimes \sO_Y(-(\nu-1)\cdot H)
$$
over $U_0$. By (\ref{2.7.3}) one finds morphisms, surjective over $U_0$
\begin{multline*}
f^{(r)}_*(\delta^*(\omega_{X/Y}^{\nu-1}\otimes
\omega_{X^{(r)}/Y}))\otimes \sO_Y(-(\nu-1)\cdot H)\\
\>>> f^{r}_*(\omega_{X^r/Y}^{\nu})\otimes \sO_Y(-(\nu-1)\cdot H)
=\Big( \bigotimes^rf_*\omega_{X/Y}^{\nu}\Big)
\otimes \sO_Y(-(\nu-1)\cdot H)\\
\>>>
S^r(f_*\omega_{X/Y}^{\nu})
\otimes \sO_Y(-(\nu-1)\cdot H).
\end{multline*}
Since the quotient of a weakly positive sheaf is weakly
positive, the sheaf on the right hand side is weakly positive
over $U_0$, hence $f_*\omega_{X/Y}^{\nu}$ is ample with respect
to $U_0$. For $\nu$ sufficiently large $\sO_Y((\nu-1)\cdot H -S)$
is ample, and one obtains the second part of \ref{2.6}.
\end{proof}

If $f:X\to Y$ is not semi-stable in codimension one, the
sheaf of relative $n$-forms $\Omega^n_{X/Y}(\log \Delta)$ might
be strictly smaller than the relative dualizing sheaf
$\omega_{X/Y}$. In fact, comparing the first Chern classes of
the entries in the tautological sequence
\begin{equation}\label{2.7.6}
0 \>>> f^*\Omega^1_Y(\log S) \>>> \Omega^1_X(\log \Delta)
\>>> \Omega^1_{X/Y}(\log \Delta) \>>> 0
\end{equation}
one finds for $\Delta=f^* S$
\begin{equation}\label{2.7.7}
\Omega^n_{X/Y}(\log \Delta) = \omega_{X/Y}(\Delta_{\rm red} - \Delta).
\end{equation}

\begin{corollary} \label{2.8}
Under the assumptions made in \ref{2.6}, i), ii) or iii), for all
$\nu$ sufficiently large and divisible, the sheaf $f_*\Omega^n_{X/Y}(\log
\Delta)^\nu$ is ample with respect to $U_0$.
\end{corollary}

Before proving \ref{2.8} let us start to study
the behavior of the relative $q$-forms under base extensions.
Here we will prove a more general result than needed for
\ref{2.8}, and we will not require $Y\setminus U$ to be smooth.
\begin{assumptions}\label{2.9a}
Let $f:X\to Y$ be any partial compactification of a smooth
family $V\to U$, let $\psi:Y'\to Y$ be a finite covering with $Y'$ non
singular, and let $\tilde{X}$ be the normalization of
$X\times_YY'$. Consider a desingularization $\varphi:X' \to
\tilde{X}$, where we assume the center of $\varphi$ to lie in
the singular locus of $\tilde{X}$.
The induced morphisms are denoted by
\begin{equation}\label{2.9.1}
\begin{CDS}
X' \> \varphi >> \tilde{X} \> \tilde{\varphi} >> \hspace{-.2cm}X\times_YY'
\hspace{-.2cm}\> \pi_1 >> X \\
\novarr \SE f'\ \ E E \V \tilde{f} V V \SW \pi_2 W W \novarr \SW W f W \\
\noharr Y' \> \psi >> \ \ Y.
\end{CDS}
\end{equation}
Let us define $\delta=\tilde{\varphi}\circ \varphi: X' \to X\times_YY'$ and
$\psi'=\pi_1\circ\delta:X'\to X$. Finally we
write $S'=\psi^* S$ and $\Delta'={\psi'}^*\Delta$.
The discriminant loci of $\psi$ and $\psi'$ will be
$\Delta(X'/X)$, and $\Delta(Y'/Y)$, respectively. We will assume that
$ S+ \Delta(Y'/Y)$ and $\Delta+\Delta(X'/X)$, as well as their
preimages in $Y'$ and $X'$, are normal crossing divisors.
\end{assumptions}
\begin{lemma}\label{2.9}
Using the assumptions and notations from \ref{2.9a},
\begin{enumerate}
\item[i)] there exists for all $p$ an injection
$$
{\psi'}^*\Omega^p_{X/Y}(\log \Delta)\> \subset >>
\Omega^p_{X'/Y'}(\log \Delta'),
$$
which is an isomorphism over ${\psi'}^{-1}(X\setminus {\rm
Sing}(\Delta))$.
\item[ii)] there exists for all $\nu >0$ an injection
\begin{gather*}
f'_*(\Omega^n_{X'/Y'}(\log \Delta')^{(\nu-1)}\otimes
\omega_{X'/Y'}) \>>> \psi^*f_*(\Omega^n_{X/Y}(\log
\Delta)^{(\nu-1)}\otimes\omega_{X/Y}) ,
\end{gather*}
which is an isomorphism over $\psi^{-1}(U)$.
\end{enumerate}
\end{lemma}
\begin{proof}
If one replaces in the tautological sequence (\ref{2.7.6})
the divisor $S$ by a larger one, the sheaf on the right hand side does not
change, hence
$$
\Omega^1_{X/Y}(\log \Delta) = \Omega^1_{X/Y}(\log
(\Delta+\Delta(X'/X))).
$$
Both, $\Omega^1_Y(\log (S + \Delta(Y'/Y)))$ and $\Omega^1_{X}(\log
(\Delta+\Delta(X'/X)))$ behave well under pullback to $X'$ (see
\cite{E-V}, 3.20, for example). To be more precise,
there exists an isomorphism
$$
\psi^*\Omega^1_Y(\log (S + \Delta(Y'/Y)))\simeq
\Omega^1_{Y'}(\log (S' + \psi^*\Delta(Y'/Y)))
$$
and an injection
$$
{\psi'}^*\Omega^1_{X}(\log (\Delta+\Delta(X'/X)))\> \subset >>
\Omega^1_{X'}(\log {\psi'}^*(\Delta+\Delta(X'/X)))
$$
which is an isomorphism over the largest open subscheme $V_1'$, where
$\psi'$ is an isomorphism. Since $\tilde{X}$ is non-singular
outside of $\Delta$, and since the singularities of $\tilde{X}$
can only appear over singular points of the discriminant
$\Delta(X'/X)$, we find ${\psi'}^{-1}(X\setminus {\rm
Sing}(\Delta))\subset V_1'$.

For ii) we use again that $\tilde{X}$ is non-singular over $X\setminus {\rm
Sing}(\Delta)$. So part i) induces an isomorphism
$$
\varphi_*(\Omega^n_{X'/Y'}(\log \Delta')^{(\nu-1)}\otimes
\omega_{X'/Y'}) \>\simeq >> \tilde{\varphi}^*\pi_1^*(\Omega^n_{X/Y}(\log
\Delta)^{(\nu-1)})\otimes\omega_{\tilde{X}/Y'}.
$$
The natural map $\tilde{\varphi}_*\omega_{\tilde{X}/Y'} \>>>
\omega_{X\times_YY'/Y'}$ and the projection formula
give
$$
\delta_*(\Omega^n_{X'/Y'}(\log \Delta')^{(\nu-1)}\otimes
\omega_{X'/Y'}) \>>> \pi_1^*(\Omega^n_{X/Y}(\log
\Delta)^{(\nu-1)}\otimes\omega_{X/Y}),
$$
and ii) follows by flat base change.
\end{proof}
\noindent
{\it Proof of \ref{2.8}.}
Applying \ref{1.2}, b), one finds a finite covering $\psi:Y'\to Y$
such that the family $f':X'\to Y'$ is semi-stable in
codimension one, hence (\ref{2.7.7}) implies
$\omega_{X'/Y'}=\Omega^n_{X'/Y'}(\log \Delta')$,
whereas $\omega_{X/Y} \otimes f^*\sO_Y(-S)\subset
\Omega^n_{X/Y}(\log (\Delta))$.
So \ref{2.9}, ii), gives a morphism of sheaves
\begin{multline*}
f'_*(\omega_{X'/Y'}^\nu)\otimes \sO_{Y'}(-\psi^*S) \>>>
\psi^*(f_*(\Omega^n_{X/Y}(\log
\Delta)^{(\nu-1)}\otimes\omega_{X/Y})\otimes \sO_Y(-S))\\
\>>> \psi^*f_*(\Omega^n_{X/Y}(\log \Delta)^\nu)
\end{multline*}
By \ref{2.6}, iii), for some $\nu \gg 0$ the sheaf on the left
hand side will be ample with respect to $\psi^{-1}(U_0)$, hence the
sheaf on the right hand side has the same property.
\qed\\
\ \\
\noindent
A positivity property, similar to the last one, will be
expressed in terms of fibred products of the given family. It will
be used in the proof of \ref{0.3}, iii). We do not need it in its
full strength, just ``up to codimension two in $Y$''.
Nevertheless, in order to be able to refer to \cite{V-Z2}, we
formulate it in a more general setup.

Let $V\to U$ be a smooth family with $\omega_{V/U}$
$f$-semi-ample. By \ref{1.5} we find a smooth birational model
$\tilde{V}\to\tilde{U}$ whose compactification $f:X\to Y$
fits into the diagram (\ref{1.5.1}):
$$
\begin{CDS}
X \< \psi' << X' \< \sigma << Z \\
\V V f V \novarr \V V f' V \SW W g W \\
Y \< \psi << Y'\novarr
\end{CDS}
$$
Let us choose any $\nu \geq 3$ such that
$$
f^*f_*\omega_{X/Y}^\nu \>>> \omega_{X/Y}^\nu
$$
is surjective over $\tilde{V}$,
and that the multiplication map
$$
S^\eta(f_*\omega^\nu_{X/Y}) \>>> f_*\omega^{\eta\cdot \nu}_{X/Y}
$$
is surjective over $\tilde{U}$. By definition one has ${\rm
Var}(f)={\rm Var}(g)$. If ${\rm Var}(f)=\dim(Y)$,
applying \ref{2.6}, i), to $g$ one finds that the sheaf
$\lambda_\nu$, defined in \ref{1.5}, d), is of maximal Kodaira
dimension. Hence some power of $\lambda_\nu$ is of
the form $\sA(D)$, for an ample invertible sheaf $\sA$ on $Y$ and
for an effective divisor $D$ on $Y$. We may assume moreover,
that $D \geq S$ and, replacing the number $N_\nu$ in \ref{1.5}
by some multiple, that
$$
\det (g_* \omega^{\nu}_{Z/Y'})^{N_{\nu}} =
\sA(D)^{\nu\cdot (\nu-1)\cdot e}
$$
where $e = {\rm Max}\{e(\omega_F^\nu);$ $F$ a fibre of $V\to U\}$.

\begin{proposition} \label{2.10}
For $r=N_\nu\cdot {\rm rank}(f_*\omega_{X/Y}^\nu)$,
let $X^{(r)}$ denote a desingularization of the $r$-th fibre
product $X \times_Y \ldots \times_Y X$ and let $f^{(r)} : X^{(r)}
\to Y$ be the induced family. Then for all $\beta$
sufficiently large and divisible the sheaf
$$
f^{(r)}_*(\Omega^{r\cdot n}_{X^{(r)}/Y}(\log \Delta)^{\beta\cdot\nu}) \otimes
\sA^{-\beta\cdot\nu\cdot (\nu-2)}
$$
is globally generated over some non-empty open subset $U_0$  of
$\tilde{U}$, and the sheaf
$$
\Omega^{r\cdot n}_{X^{(r)}/Y}(\log \Delta)^{\beta\cdot\nu} \otimes
f^{(r)*}\sA^{-\beta\cdot\nu\cdot (\nu-2)}
$$
is globally generated over ${f^{(r)}}^{-1}(U_0)$.
\end{proposition}
\begin{proof}
By \cite{V-Z2}, 4.1, the sheaf
$$
f^{(r)}_*(\omega_{X^{(r)}/Y}^{\beta\cdot\nu}) \otimes
\sA^{-\beta\cdot\nu\cdot (\nu-2)}\otimes
\sO_Y(-\beta\cdot\nu\cdot(\nu-1)\cdot D)
$$
is globally generated over some open subset.
However, by (\ref{2.7.7})
$$
\omega_{X^{(r)}/Y}^{\beta\cdot\nu}\otimes
f^*\sO_Y(-\beta\cdot\nu \cdot S)
$$
is contained in
$$
\Omega^{r\cdot n}_{X^{(r)}/Y}(\log \Delta)^{\beta\cdot\nu}.
$$
Since
$$
\beta\cdot\nu\cdot(\nu-1)\cdot D \geq \beta\cdot\nu \cdot S
$$
one obtains \ref{2.10}, as stated.
\end{proof}

\section{Higgs bundles and the proof of \ref{0.3}}
\label{3}

As in \cite{V-Z} and \cite{V-Z2}, in order to prove \ref{0.3}
we have to construct certain Higgs bundles, and we have to
compare them to one, induced by a variation of Hodge structures.
For \ref{0.3}, iii), we will just use the content of the second
half of section \ref{2}. For iv) we need in addition Kawamata's
covering construction, as explained in \ref{1.2}. The
reduction steps contained in the second half of
section \ref{1} will be needed for \ref{0.3}, i).

So let $U$ be a manifold and let $Y$ be a smooth projective
compactification with $Y\setminus U$ a normal crossing divisor.
Starting with a smooth family $V\to U$ with $\omega_{V/U}$
relative semi-ample over $U$, we first choose a smooth
projective compactification $X$ of $V$, such that $V\to U$ extends
to $f:X\to Y$.

In the first half of the section, we will work with
good partial compactifications as
defined in \ref{1.1}. Hence leaving out a codimension
two subscheme of $Y$, we will assume that the
divisor $S=Y\setminus U$ is smooth, that $f$ is flat and
that $\Delta=X\setminus V$ is a relative normal crossing divisor.
The exact sequence (\ref{2.7.6}) induces a filtration on
the wedge product $\Omega^p_{X/Y}(\log \Delta)$, and thereby
the tautological sequences
\begin{multline}\label{3.0.1}
0\to {f}^*\Omega^1_Y(\log S)\otimes
\Omega^{p-1}_{X/Y}(\log \Delta)
\to {\mathfrak g \mathfrak r}(\Omega_X^p(\log
\Delta)) \to \Omega_{X/Y}^p(\log \Delta)\to 0,
\end{multline}
where
\begin{gather*}
{\mathfrak g \mathfrak r}(\Omega_X^p(\log \Delta))=
\Omega_X^p(\log \Delta)
/f^*\Omega^2_Y(\log S)\otimes \Omega^{p-2}_{X/Y}(\log
\Delta).
\end{gather*}
Given an invertible sheaf $\sL$ on $X$ we will study in this
section various sheaves of the form
$$
F_0^{p,q}:=R^qf_*(\Omega^{p}_{X/Y}(\log
\Delta)\otimes\sL^{-1})/_{\rm torsion}
$$
together with the edge morphisms
$$
\tau^0_{p,q}:F_0^{p,q}\>>> F_0^{p-1,q+1}\otimes \Omega^1_{Y}(\log
S),
$$
induced by the exact sequence (\ref{3.0.1}), tensored with
$\sL^{-1}$.

First we have to extend the base change properties
for direct images, studied in \ref{2.9}, ii), to higher
direct images.
\begin{lemma}\label{3.1}
Keeping the notations and assumptions from \ref{2.9a}, let
$Y'_1$ be the largest open subset in $Y'$ with $X\times_YY'_1$ normal.
We write
$$
\iota: \psi^*\Omega_Y^1(\log S) \>>> \Omega_{Y'}^1(\log S')
$$
for natural inclusion, and we consider an invertible sheaf $\sL$
on $X$, and its pullback $\sL'={\psi'}^*\sL$ to $X'$.

Then for all $p$ and $q$, there are morphisms
$$
{\psi}^*F_0^{p,q} \> \zeta_{p,q} >>
{F'_0}^{p,q}:=R^qf'_*(\Omega^{p}_{X'/Y'}(\log
\Delta')\otimes{\sL'}^{-1})/_{\rm torsion},
$$
whose restriction to $Y'_1$ are isomorphisms,
and for which the diagram
\begin{equation}\label{3.1.1}
\begin{CD}
{\psi}^*F_0^{p,q} \> \psi^*(\tau^0_{p,q}) >>
{\psi}^*F_0^{p-1,q+1}\otimes\Omega_Y^1(\log S))\\
\V \zeta_{p,q} VV \V \zeta_{p-1,q+1}\otimes\iota VV\\
{F'_0}^{p,q} \>{\tau'}^0_{p,q} >>
{F'_0}^{p-1,q+1}\otimes\Omega_{Y'}^1(\log S')
\end{CD}
\end{equation}
commutes. Here ${\tau'}^0_{p,q}$ is again the edge morphism
induced by the exact sequence on $X'$, corresponding to
(\ref{3.0.1}) and tensored with ${\sL'}^{-1}$.
\end{lemma}
\begin{proof}
We use the notations from (\ref{2.9.1}), i.e.
\begin{equation*}
\begin{CDS}
X' \> \varphi >> \tilde{X} \> \tilde{\varphi} >> \hspace{-.2cm}X\times_YY'
\hspace{-.2cm}\> \pi_1 >> X \\
\novarr \SE f'\ \ E E \V \tilde{f} V V \SW \pi_2 W W \novarr \SW W f W \\
\noharr Y' \> \psi >> \ \ Y
\end{CDS}
\end{equation*}
and $\psi'=\varphi\circ\tilde{\varphi}\circ\pi_1$.
As in the proof of \ref{2.9}, in order to show the existence of
the morphisms $\zeta_{p,q}$ and the commutativity of the diagram
(\ref{3.1.1}) we may enlarge $S$ and $S'$ to include the discriminant loci,
hence assume that
$$
\psi^*\Omega^1_Y(\log S) = \Omega^1_{Y'}(\log S').
$$
By the generalized Hurwitz formula \cite{E-V}, 3.21,
$$
{\psi'}^*\Omega^p_X(\log \Delta) \subset \Omega^p_{X'}(\log \Delta'),
$$
and by \cite{E-V1}, Lemme 1.2,
$$
R^q\varphi_*\Omega^{p}_{X'}(\log \Delta') =
\left\{ \begin{array}{ll}
\tilde{\varphi}^*\pi_1^*\Omega^{p}_{X}(\log \Delta) & \mbox{for \ } q=0 \\
0 & \mbox{for \ } q> 0.
\end{array} \right.
$$
The tautological sequence
\begin{multline*}
0\to {f'}^*\Omega^1_{Y'}(\log (S'))\to
\Omega_{X'}^1(\log \Delta') \to
\Omega_{X'/Y'}^1(\log \Delta') \to 0
\end{multline*}
defines a filtration on $\Omega_{X'}^p(\log \Delta')$,
with subsequent quotients isomorphic to
$$
{f'}^*\Omega_{Y'}^\ell(\log S') \otimes
\Omega_{X'/Y'}^{p-\ell}(\log \Delta').
$$
Induction on $p$ allows to deduce that
\begin{equation}\label{3.1.2}
R^q\varphi_*\Omega^{p}_{X'/Y'}(\log \Delta') =
\left\{ \begin{array}{ll}
\tilde{\varphi}^*\pi_1^*\Omega^{p}_{X/Y}(\log \Delta) & \mbox{for \ } q=0 \\
0 & \mbox{for \ } q> 0.
\end{array} \right.
\end{equation}
On the other hand, the inclusion $\sO_{Z\times_YY'} \to
\tilde{\varphi}_*\sO_{\tilde{Z}}$ and flat base change gives
\begin{multline}\label{3.1.3} \hspace{1cm}
\psi^*F_0^{p,q}=
{\psi}^*R^qf_*(\Omega^{p}_{X/Y}(\log \Delta)\otimes\sL^{-1})
\> \simeq >>\\ R^q{\pi_2}_*(\pi_1^*(\Omega^{p}_{X/Y}(\log
\Delta)\otimes\sL^{-1}))
\to R^q{\tilde{f}}_*(\tilde{\varphi}^*(\pi_1^*(\Omega^{p}_{X/Y}(\log
\Delta)\otimes\sL^{-1}))) = {F'_0}^{p,q},
\end{multline}
hence $\zeta_{p,q}$. The second morphism in (\ref{3.1.3}) is an
isomorphism on the largest open subset where $\tilde{\varphi}$
is an isomorphism, in particular on $\tilde{f}^{-1}(Y'_1)$.

The way we obtained (\ref{3.1.2}) the morphisms are obviously
compatible with the different tautological sequences.
Since we assumed $S$ to contain the discriminant
locus, the pull back of (\ref{3.0.1}) to $\tilde{X}$ is
isomorphic to
\begin{multline*}
0\to \varphi_*({f'}^*\Omega^1_{Y'}(\log S')\otimes
\Omega^{p-1}_{X'/Y'}(\log \Delta'))
\to \varphi_*({\mathfrak g \mathfrak r}(\Omega_{X'}^p(\log
\Delta')))\\
\to \varphi_*(\Omega_{X'/Y'}^p(\log \Delta'))\to 0,
\end{multline*}
and the diagram (\ref{3.1.1}) commutes.
\end{proof}
\begin{remark}\label{3.2}
If $\psi:Y'\to Y$ is any smooth morphism, then again
$X\times_YY'$ non-singular.
The compatibility of the $F_0^{p,q}$ with pullback, i.e.
the existence of an isomorphism $\zeta_{p,q}:\psi^*F_0^{p,q} \to
{F'_0}^{p,q}$, and the commutativity of (\ref{3.1.1}) is also
guaranteed, in this case. In fact, both $\varphi$ and $\tilde{\varphi}$
are isomorphisms, as well as the two morphisms in
(\ref{3.1.3}).
\end{remark}
\begin{corollary}\label{3.3}
Keeping the assumptions from \ref{3.1}, assume that
$X\times_YY'$ is normal. Then the image of
$$
{F'_0}^{p,q} \>{\tau'}^0_{p,q} >>
{F'_0}^{p-1,q+1}\otimes\Omega_{Y'}^1(\log S')
$$
lies in ${F'_0}^{p-1,q+1}\otimes\psi^*(\Omega_{Y}^1(\log S))$.
\end{corollary}
\qed\\

In the sequel we will choose $\sL=\Omega^n_{X/Y}(\log \Delta).$
Let us consider first the case that for some $\nu \gg 1$ and for
some invertible sheaf $\sA$ on $Y$ the sheaf
\begin{gather}
\sL^\nu\otimes f^*\sA^{-\nu} \mbox{ \ \ is globally generated
over \ \ } V_0=f^{-1}(U_0), \label{3.3.1}
\end{gather}
for some open dense subset $U_0$ of $Y$.

We will recall some of the constructions performed
in \cite{V-Z2}, \S 6.
Let $H$ denote the zero divisor of a section of $\sL^\nu\otimes
f^*\sA^{-\nu}$, whose restriction to a general fibre of $f$ is
non-singular. Let $T$ denote the closure of the discriminant of
$H\cap V \to U$. Leaving out some more codimension two
subschemes, we may assume that $S+T$ is a smooth divisor. We
will write $\Sigma = f^*T$ and we keep the notation
$\Delta=f^*(S)$.

Let $\delta:W\to X$ be a blowing up of $X$ with centers in
$\Delta+\Sigma$ such that $\delta^*(H+\Delta+\Sigma)$ is a
normal crossing divisor. We write
$$
\sM=\delta^*(\Omega^n_{X/Y}(\log \Delta)\otimes
f^*\sA^{-1}).
$$
Then for $B=\delta^*H$ one has $\sM^\nu=\sO_W(B)$. As in
\cite{E-V}, \S 3, one obtains a cyclic covering of $W$, by taking
the $\nu$-th root out of $B$. We choose $Z$ to be a
desingularization of this covering
and we denote the induced morphisms by
$g:Z\to Y$, and  $h:W\to Y$.
Writing $\Pi=g^{-1}(S\cup T)$, the restriction of $g$ to $Z_0=
Z\setminus\Pi$ will be smooth.

For the normal crossing divisor $B$ we define
$$
{\sM}^{(-1)}={\sM}^{-1}\otimes
\sO_{W}\Big(\Big[\frac{B}{\nu}\Big]\Big), \mbox{ \ \ and \ \ }
{\sL}^{(-1)}= \delta^*({\sL}^{-1})\otimes
\sO_{W}\Big(\Big[\frac{B}{\nu}\Big]\Big).
$$
In particular the cokernel of the inclusion $\delta^*\sL^{-1}\subset
\sL^{(-1)}$ lies in $h^{-1}(S+T)$. The sheaf
\begin{multline*}
F^{p,q} = R^qh_*(\delta^*(\Omega_{X/Y}^p(\log \Delta))\otimes
{\sM}^{(-1)}) \otimes\sA^{-1}/_{\rm torsion}\\
= R^qh_*(\delta^*(\Omega_{X/Y}^p(\log \Delta))\otimes
{\sL}^{(-1)})/_{\rm torsion}
\end{multline*}
contains the sheaf $F_0^{p,q}$ and both are isomorphic outside
of $S+T$. The edge morphism
\begin{gather*}
\tau_{p,q}: F^{p,q} \>>> F^{p-1,q+1}\otimes \Omega^1_Y(\log S)
\end{gather*}
given by the tautological exact sequence
\begin{multline*}
0\to {h}^*\Omega^1_Y(\log S)\otimes
\delta^*(\Omega^{p-1}_{X/Y}(\log \Delta))\otimes {\sL}^{(-1)}
\to \\
\delta^*({\mathfrak g \mathfrak r}(\Omega_X^p(\log
\Delta)))\otimes {\sL}^{(-1)} \to
\delta^*(\Omega_{X/Y}^p(\log \Delta))\otimes {\sL}^{(-1)} \to 0
\end{multline*}
is compatible with $\tau^0_{p,q}$. Let us remark, that the
sheaves $F^{p,q}$ depend on the choice of the divisor $H$ and
they can only be defined assuming (\ref{3.3.1}).

Up to now, we constructed two Higgs bundles
$$
F_0=\bigoplus F_0^{p,q} \> \subset >> F=\bigoplus F^{p,q}.
$$
We will see below, that $\sA\otimes F$ can be compared with a
Higgs bundle $E$, given by a variation of Hodge structures.
This will allow to use the negativity of the kernel of
Kodaira-Spencer maps (see \cite{Zuo}), to show that ${\rm
Ker}(\tau_{p,q})^\vee$ is big.

By \cite{Del}, for all $k\geq 0$, the local constant system $R^kg_*\C_{Z_0}$
gives rise to a local free sheaf $\sV_k$ on $Y$ with the
Gau\ss-Manin connection
$$
\nabla: \sV_k \>>> \sV_k \otimes\Omega^1_Y(\log(S+T)).
$$
We assume that $\sV_k$ is the quasi-canonical extension of
$$
(R^kg_*\C_{Z_0})\otimes_\C \sO_{Y\setminus (S\cup T)},
$$
i.e. that the real part of the eigenvalues of the
residues around the components of $S+T$ lie in $[0,1)$.

Since we assumed $S+T$ to be non-singular, $\sV_k$ carries a
filtration $\sF^p$ by subbundles (see \cite{Sch}).
So the induced graded sheaves $E^{p,k-p}$ are locally free,
and they carry a Higgs structure with logarithmic poles along
$S+T$. Let us denote it by
$$
({\mathfrak g \mathfrak r}_{\sF}(\sV_k), {\mathfrak g \mathfrak
r}_{\sF}(\nabla))=(E,\theta)=\Big(
\bigoplus_{q=0}^{k}E^{k-q,q}\ , \
\bigoplus_{q=0}^{k}\theta_{k-q,q}\Big).
$$
As well-known (see for example \cite{Gri}, page 130) the bundles
$E^{p,q}$ are given by
$$
E^{p,q}=R^qg_*\Omega_{Z/Y}^p(\log \Pi).
$$
Writing again ${\mathfrak g \mathfrak r}(\underline{
\ \ })$ for ``modulo the pullback of $2$-forms on $Y$'',
the Gau\ss-Manin connection is the edge morphism
of
\begin{multline*}
0\to {g}^*\Omega^1_Y(\log (S+T))\otimes
\Omega^{\bullet-1}_{Z/Y}(\log \Pi)\to
{\mathfrak g \mathfrak r}(\Omega_Z^{\bullet}(\log \Pi)) \to
\Omega_{Z/Y}^{\bullet}(\log \Pi) \to 0.
\end{multline*}
Hence the Higgs maps
\begin{gather*}
{\theta}_{p,q}:E^{p,q} \>>> E^{p-1,q+1}\otimes \Omega^1_Y(\log (S+T))
\end{gather*}
are the edge morphisms of the tautological exact sequences
\begin{multline*}
0\to {g}^*\Omega^1_Y(\log (S+T))\otimes
\Omega^{p-1}_{Z/Y}(\log \Pi)\to\\
{\mathfrak g \mathfrak r}(\Omega_Z^p(\log \Pi)) \to
\Omega_{Z/Y}^p(\log \Pi) \to 0.
\end{multline*}
In the sequel we will write $T_*(-\log **)$ for the dual of
$\Omega^1_*(\log **)$.
\begin{lemma}\label{3.4} Under the assumption (\ref{3.3.1}) and
using the notations introduced above,
let
$$
\iota:\Omega^1_Y(\log S)\>>>\Omega^1_Y(\log (S+T))
$$
be the natural inclusion. Then there exist morphisms
$\rho_{p,q}: \sA\otimes F^{p,q} \to E^{p,q}$
such that:
\begin{enumerate}
\item[i)] The diagram
$$
\begin{CD}
E^{p,q} \> \theta_{p,q} >> E^{p-1, q+1} \otimes \Omega^{1}_{Y}
(\log (S+T)) \\
\A \rho_{p,q} AA \A A \rho_{p-1,q+1} \otimes \iota A \\
\sA\otimes F^{p,q} \> {\rm id}_{\sA}\otimes\tau_{p,q}
>> \sA\otimes F^{p-1,q+1} \otimes \Omega^{1}_{Y} (\log S).
\end{CD}
$$
commutes.
\item[ii)] $F^{n,0}$ has a section $\sO_Y \to F^{n,0}$, which is
an isomorphism on $Y\setminus (S\cup T)$.
\item[iii)] $\tau_{n,0}$ induces a morphism
$$
\tau^\vee:T_Y( - \log S)=(\Omega_Y^1(\log S))^\vee \>>>
F^{{n,0}^\vee}\otimes F^{n-1,1},
$$
which coincides over $Y\setminus (S\cup T)$ with the Kodaira-Spencer map
$$
T_Y( - \log S) \>>> R^1f_*T_{X/Y}(-\log \Delta).
$$
\item[iv)] $\rho_{n,0}$ is injective.
If the general fibre of $f$ is canonically polarized,
then the morphisms $\rho_{n-m,m}$ are injective, for all $m$.
\item[v)] Let $\sK^{p,q}={\rm Ker}(E^{p,q} \> \theta_{p,q} >>
E^{p-1, q+1} \otimes \Omega^{1}_{Y} (\log (S+T)))$. Then
the dual $(\sK^{p,q})^\vee$ is weakly positive with respect to
some open dense subset of $Y$.
\item[vi)] The composite
\begin{equation*}
\theta_{n-q+1,q-1}\circ \cdots \circ \theta_{n,0}:
E^{n,0} \>>> E^{n-q,q}\otimes \bigotimes^{q} \Omega^1_Y(\log
(S+T))
\end{equation*}
factors like
\begin{equation*} \ \ \ \ \ \ \
E^{n,0} \> \theta^q >> E^{n-q,q}\otimes S^{q} \Omega^1_Y(\log
(S+T)) \> \subset >> E^{n-q,q}\otimes \bigotimes^{q} \Omega^1_Y(\log
(S+T)).
\end{equation*}
\end{enumerate}
\end{lemma}
\begin{proof}
The properties i) - iv) have been verified in \cite{V-Z2}, 6.3
in case the general fibre is canonically polarized. So let us
just sketch the arguments.

By \cite{E-V} (see also \cite{V-Z2}, 6.2) the sheaf
$$
R^q{h}_*(\Omega_{W/Y}^p(\log
(B+\delta^*\Delta+\delta^*\Sigma))\otimes {\sM}^{(-1)})
$$
is a direct factor of $E^{p,q}$. The morphism $\rho_{p,q}$
is induced by the natural inclusions
\begin{multline*}
\delta^*\Omega_{X/Y}^p(\log \Delta)\to
\delta^*\Omega_{X/Y}^p(\log (\Delta+\Sigma))\\
\to \Omega_{W/Y}^p(\log (\delta^*\Delta+\delta^*\Sigma))\to
\Omega_{W/Y}^p(\log (B+\delta^*\Delta+\delta^*\Sigma)),
\end{multline*}
tensored with $\sM^{(-1)}=\sL^{(-1)}\otimes h^* \sA$.

Such an injection also exist for $Y$ replaced by ${\rm
Spec}(\C)$. Since the different tautological sequences are compatible
with those inclusions one obtains i).
Over $Y\setminus (S\cup T)$ the kernel of $\rho_{n-m,m}$ is a quotient of
the sheaf
$$
R^{m-1}(h|_{B})_*(\Omega_{B/Y}^{n-m-1} \otimes {\sM}^{-1}|_{B}).
$$
In particular $\rho_{n,0}$ is injective. The same holds true for
all the $\rho_{n-m,m}$ in case $\sM$ is fibre wise ample,
by the Akizuki-Kodaira-Nakano vanishing theorem.

By definition
$$
F^{n,0}= h_*(\delta^*(\Omega^n_{X/Y}(\log \Delta))\otimes
{\sL}^{(-1)}) = h_*\sO_{W}\Big(\Big[\frac{B}{\nu}\Big]\Big),
$$
and ii) holds true.

For iii), recall that over $Y\setminus (S\cup T)$
the morphism
$$
\delta^*(\sL \otimes f^* \sA^{-1})=\sM^{-1} \to {\sM}^{(-1)}
$$
is an isomorphism. By the projection formula the morphism
$\tau_{n,0}|_{Y\setminus (S\cup T)}$ is the restriction of the
edge morphism of the short exact sequence
\begin{equation*}
0\to {f}^*\Omega^1_U \otimes \Omega^{n-1}_{V/U}\otimes \sL^{-1}\to
{\mathfrak g \mathfrak   r}(\Omega_{V}^n)\otimes \sL^{-1} \to
\Omega_{V/U}^{n} \otimes \sL^{-1} \to 0.
\end{equation*}
The sheaf on the right hand side is $\sO_V$ and the one on the left hand side
is $f^*\Omega^1_U\otimes T_{V/U}$.
For $r=\dim(U)$, tensoring the exact sequence with
$$
f^*T_U=f^*(\Omega_U^{r-1} \otimes \omega_U^{-1})
$$
and dividing by the kernel of the wedge product
$$
f^*\Omega_U^1\otimes f^*(\Omega_U^{r-1} \otimes \omega_U^{-1})
=f^*\Omega_U^1\otimes f^*T_U \>>> \sO_V
$$
on the left hand side, one obtains an exact
sequence
\begin{equation}\label{3.4.1}
0 \>>> T_{V/U} \>>> \sG \>>> f^*T_U  \>>> 0,
\end{equation}
where $\sG$ is a quotient of
${\mathfrak g \mathfrak r}(\Omega_{V}^n)\otimes
{\omega_{V}}^{-1} \otimes f^* \Omega^{r-1}_U$.
By definition, the restriction to $Y\setminus (S\cup T)$ of
the morphism considered in iii) is the first edge
morphism in the long exact sequence, obtained by applying
$R^\bullet f_*$ to (\ref{3.4.1}).

The wedge product induces a morphism
$$
\Omega_{V}^n \otimes {\omega_{V}}^{-1} \otimes f^*
\Omega^{r-1}_U \>>> \Omega_V^{n+r-1} \otimes \omega_V^{-1}=T_V.
$$
This morphism factors through $\sG$. Hence the exact sequence
(\ref{3.4.1}) is isomorphic to the tautological sequence
\begin{equation}\label{3.4.2}
0\>>> T_{V/U} \>>> T_V \>>> f^*T_U \>>> 0.
\end{equation}
The edge morphism $T_U \to R^1f_*T_{V/U}$ of (\ref{3.4.2})
is the Kodaira-Spencer map.

In order to prove v), we use as in the proof of \ref{1.2}, b),
Kawamata's covering construction to find a non-singular finite
covering $\rho: Y' \to Y$ such that for some desingularization
$Z'$ of $Z\times_YY'$ the induced variation of Hodge structures
has uni-potent monodromy, and such that $g':Z'\to Y'$ is semi-stable.

From \ref{3.1}, applied to $Z$, $\sO_Z$ and $\Pi$ instead of
$X$, $\sL$ and $\Delta$ one obtains a commutative diagram
$$
\begin{CD}
\rho^*{E}^{p,q} \> \rho^*\theta_{p,q}>> \rho^*{E}^{p-1,q+1}\otimes
\Omega^1_{Y} (\log S+T)\\
\V\subset VV \V\subset VV\\
{E'}^{p,q} \> \theta'_{p,q}>> {E'}^{p-1,q+1}\otimes
\Omega^1_{Y'} (\log S'),
\end{CD}
$$
where $S'=\psi^* (S+T)$, where ${E'}^{p,q} =
R^qg'_*\Omega^p_{Z'/Y'}(\log \Pi'),$ and where $\theta'_{p,q}$
is the edge-morphism.

In particular the pullback of the kernel of
$\theta_{p,q}$, the sheaf $\rho^*\sK^{p,q}$, lies in the kernel
${\sK'}^{p,q}$ of $\theta'_{p,q}$. Leaving out some codimension
two subschemes of $Y$ and $Y'$, we may assume that ${\sK'}^{p,q}$
is a subbundle of ${E'}^{p,q}$. Choose a smooth extension
$\bar{Y'}$ of $Y'$ such that the closure of
$S'\cup(\bar{Y'}-Y')$ is a normal crossing
divisor, and let ${\bar{E'}}^{p,q}$ be the Higgs bundle,
corresponding to the canonical extension of the variation of
Hodge structures. For some choice of the compactification
${\sK'}^{p,q}$ will extend to a subbundle ${\bar{\sK'}}^{p,q}$ of
${\bar{E'}}^{p,q}$. By \cite{Zuo}, 1.2, the dual
$({\bar{\sK'}}^{p,q})^\vee$ is
numerically effective, hence weakly positive. Thereby
$\rho^*(\sK^{p,q})^\vee$ is weakly positive over some
open subset, and the compatibility of weak positivity with
pullback shows v).

For vi) one just has to remark that on page 12 of \cite{Sim} it
is shown that $\theta\wedge\theta=0$ for
$$
\theta=\bigoplus_{q=0}^n\theta_{n-q,q}.
$$
\end{proof}
\begin{corollary}\label{3.5}
Assume (\ref{3.3.1}) holds true for some ample invertible sheaf
$\sA$ and for some $\nu \gg 1$.
Assume moreover that there exists a locally free subsheaf
$\Omega$ of $\Omega^1_Y(\log S)$ such that ${\rm
id}_{\sA}\otimes \tau_{p,q}$ factors through
$$
\sA\otimes F^{n-q,q} \>>> \sA\otimes F^{n-q-1,q+1} \otimes \Omega,
$$
for all $q$. Then for some $0< m \leq n$ there exists a big coherent subsheaf
$\sP$ of $S^m(\Omega)$.
\end{corollary}
\begin{proof}
Using the notations from \ref{3.4}, write
$\sA\otimes\tilde{F}^{n-q,q}=\rho_{n-q,q}(\sA\otimes F^{n-q,q})$.
By \ref{3.4}, i), and by the choice of $\Omega$
$$
\theta_{n-q,q}(\sA\otimes \tilde{F}^{n-q,q}) \subset
\sA\otimes \tilde{F}^{n-q-1,q+1}\otimes \Omega.
$$
By \ref{3.4}, ii) and iv), there is a section
$\sO_Y \to F^{n,0} \simeq \tilde{F}^{n,0}$, generating
$\tilde{F}^{n,0}$ over $Y\setminus (S \cup T)$,
and by \ref{3.4}, v), $\sA\otimes \tilde{F}^{n,0}$ can not lie in
the kernel of $\theta_{n,0}$. Hence the largest number $m$ with
$\theta^m (\sA\otimes\tilde{F}^{n,0})\neq 0$
satisfies $1\leq m\leq n$. By the choice of $m$
$$
\theta^{m+1} (\sA\otimes\tilde{F}^{n,0})= 0,
$$
and \ref{3.4}, vi) implies that $\theta^m
(\sA\otimes\tilde{F}^{n,0})$ lies in
$$
(\sK^{n-m,m}\cap \sA\otimes \tilde{F}^{n-m,m}) \otimes
S^m(\Omega) \subset \sK^{n-m,m}\otimes S^m(\Omega).
$$
We obtain morphisms of sheaves
$$
\sA \otimes (\sK^{n-m,m})^\vee \> \subset >> \sA\otimes
\tilde{F}^{n,0} \otimes (\sK^{n-m,m})^\vee \> \neq 0 >> S^m(\Omega).
$$
By \ref{3.4}, v), the sheaf on the left hand side is big, hence
its image $\sP\subset S^m(\Omega)$ is big as well.
\end{proof}

\noindent
{\it Proof of \ref{0.3}, iii).} Let $Y$ be the given
smooth projective compactification of $U$ with $Y\setminus U$ a
normal crossing divisor. In order to prove iii) we may blow up $Y$.
Hence given a morphism $V\to U$ with $\omega_{V/U}$ semi-ample,
by abuse of notations we will assume that $V\to U$ itself
fits into the diagram (\ref{1.5.1}).
So we may apply \ref{2.10} and replace $X$ by $X^{(r)}$ for $r$
sufficiently large. In this way we loose control on the
dimension of the fibres, but we enforce the existence of a
family for which (\ref{3.3.1}) holds true.
We obtain the big coherent subsheaf $\sP$, asked for in
\ref{0.3}, ii), by \ref{3.5}, applied to $\Omega=\Omega^1_Y(\log
S)$.
\qed\\
\ \\
In order to prove \ref{0.3}, iv), we have to argue in a slightly
different way, since we are not allowed to perform any
construction, changing the dimension of the general fibre.\\
\ \\
{\it Proof of \ref{0.3}, iv).} We start again with  a smooth
projective compactifications $X$ of $V$, such that $V\to U$
extends to $f:X\to Y$.
Recall that for $ \sL=\Omega^n_{X/Y}(\log \Delta)$, we found in \ref{2.8}
some $\nu \gg 1$ and an open dense subset $U_0$ of $Y$ such that
\begin{gather}
f_*\sL^\nu=f_*\Omega^n_{X/Y}(\log \Delta)^\nu \mbox{ \ \ is
ample with respect to \ \ } U_0\label{3.5.1}\\
\mbox{ and \ \ } f^*f_*\sL^\nu \>>> \sL^\nu \mbox{ \ \ is surjective
over \ \ } V_0=f^{-1}(U_0) \label{3.5.2}.
\end{gather}
Given a very ample sheaf $\sA$ on $Y$, lemma \ref{2.1}
implies that for some $\mu'$ the sheaf $\sA^{-1}\otimes
S^{\mu'}(f_*\sL^\nu)$ is globally generated over $U_0$.
Lemma \ref{1.2}, a), allows to find some smooth covering
$\psi:Y'\to Y$ such that $\psi^*\sA={\sA'}^{\mu}$ for an
invertible ample sheaf $\sA'$ on $Y'$ and for $\mu=\mu'\cdot\nu$.
We will show, that for this covering \ref{0.3}, iv), holds true.
To this aim, we are allowed to replace $Y$ by the complement of
a codimension two subscheme, hence assume that $f:X\to Y$
is a good partial compactification, as defined in \ref{1.1}.
In particular, we can assume $f_*\sL^\nu$ to be locally free.
Then the sheaf $\sL^{\mu}\otimes f^*\sA^{-1}$ is globally generated over
$f^{-1}(U_0)$. Let $H$ be the zero divisor
of a general section of this sheaf, and let $T$ denote the
non-smooth locus of $H \to Y$. Leaving out some additional codimension two subset,
we may assume that the discriminant $\Delta(Y'/Y)$ does not meet $T$
and the boundary divisor $S$, hence in particular that the fibred product
$X'=X\times_YY'$ is smooth. If $\psi':X'\to X$ and $f':X'\to Y'$
denote the projections, we write $S'=\psi^*S$, $T'=\psi^*T$,
$\Delta'={\psi'}^*(\Delta)={h}^*(S')$,
$$
\sL'=\Omega^n_{X'/Y'}(\log \Delta')={\psi'}^*\sL,
$$
and so on. The sheaf
$$
{\sL'}^\mu\otimes {f'}^*{\sA'}^{-\mu}= {\psi'}^*(\sL^\mu\otimes
f^*\sA^{-1})
$$
is globally generated over $\psi^{-1}(V_0)$ and (\ref{3.3.1})
holds true on $Y'$. So we can repeat the construction made
above, this time over $Y'$ and for the divisor $H'={\psi'}^*H$,
to obtain the sheaf
$$
{F'}^{p,q}= R^qh'_*{\delta'}^*(\Omega_{X'/Y'}^p(\log
\Delta')\otimes{\sL'}^{(-1)})/_{\rm torsion},
$$
together with the edge morphism
$$
\tau'_{p,q}:{F'}^{p,q}\>>>
{F'}^{p-1,q+1}\otimes\Omega^1_{Y'}(\log S'),
$$
induced by the exact sequence
\begin{multline*}
0\to {h'}^*\Omega^1_{Y'}(\log S')\otimes
{\delta'}^*(\Omega^{p-1}_{X'/Y'}(\log \Delta'))\otimes {\sL'}^{(-1)}
\to \\
{\delta'}^*({\mathfrak g \mathfrak r}(\Omega_{X'}^p(\log
\Delta')))\otimes {\sL'}^{(-1)} \to
{\delta'}^*(\Omega_{X'/Y'}^p(\log \Delta'))\otimes {\sL'}^{(-1)} \to 0.
\end{multline*}
Returning to the notations from \ref{3.1}, the sheaf
${F'_0}^{p,q}$ defined there is a subsheaf of ${F'}^{p,q}$,
both are isomorphic outside of $S'+T'$ and $\tau'_{p,q}$
commutes with ${\tau'}^0_{p,q}$. By \ref{3.3}
the image of $\tau'_{p,q}$ lies in
$\psi^*(\Omega_{Y}^1)\otimes\sO_{Y'}(*(S'+T'))$, hence in
$$
(\psi^*(\Omega^1_Y)\otimes\sO_{Y'}(*(S'+T')))\cap
\Omega_{Y'}(\log S') = \psi^*(\Omega^1_Y(\log S)).
$$
By \ref{3.5}, for some $1\leq m\leq n$ the $m$-th symmetric product of the sheaf
$$
\Omega=\psi^*(\Omega^1_Y(\log S))
$$
contains a big coherent subsheaf $\sP$, as claimed.
\qed\\
\ \\
Assume from now on, that
the fibres of the smooth family $V\to U$ are canonically
polarized, and let $f:X\to Y$ be a partial compactification.
The injectivity of $\rho_{n-m,m}$ in \ref{3.4}, iv),
gives another method to bound the number $m$ in \ref{0.3}, iii)
and to prove \ref{0.3}, ii). For i) we will use in addition, the
diagram (\ref{1.7.1}).
\begin{lemma}\label{3.6}
Using the notations from \ref{3.1}, the composite
$\tau^0_{n-q+1,q-1}\circ \cdots
\circ \tau^0_{n,0}$ factors like
$$
F_0^{n,0}=\sO_Y \> \tau_0^q >> F_0^{n-q,q}\otimes S^{q} (\Omega^1_Y(\log (S)))
\> \subset >>
F_0^{n-q,q}\otimes \bigotimes^q \Omega^1_Y(\log (S)).
$$
\end{lemma}
\begin{proof}
The equality $F_0^{n,0}=\sO_Y$ is obvious by definition. Moreover
all the sheaves in \ref{3.6} are torsion free, hence it is
sufficient to verify the existence of $\tau_0^q$ on some open
dense subset. So we may replace $Y$ by an affine subscheme, and
(\ref{3.3.1}) holds true for $\sA=\sO_Y$. By \ref{3.4}, iv), the
sheaves $F_0^{p,q}$ embed in the sheaves $E^{p,q}$, in such a way that
$\theta_{p,q}$ restricts to $\tau^0_{p,q}$. One obtains
\ref{3.6} from \ref{3.4}, vi).
\end{proof}
\noindent
{\it Proof of \ref{0.3}, i) and ii).} \ Replacing $Y$ by the
complement of a codimension two subscheme, we may choose a good
partial compactification $f:X\to Y$ of $V\to U$. Define
$$
\sN_0^{p,q} = {\rm Ker}(\tau^0_{p,q}:F_0^{p,q}\to F_0^{p-1,q+1}\otimes
\Omega_Y^1(\log S)).
$$
\begin{claim}\label{3.7}
Assume (\ref{3.3.1}) to hold true for some invertible sheaf
$\sA$, and let $(\sN_0^{p,q})^\vee$ be the dual of the sheaf
$\sN_0^{p,q}$. Then $\sA^{-1}\otimes (\sN_0^{p,q})^\vee$ is weakly
positive.
\end{claim}
\begin{proof}
Recall that under the assumption (\ref{3.3.1}) we have considered above the
slightly different sheaf
$$
F^{p,q} = R^qh_*(\delta^*(\Omega_{X/Y}^p(\log \Delta))\otimes
{\sL}^{(-1)})/_{\rm torsion},
$$
for $\delta^*\sL^{-1} \subset \sL^{(-1)}$.
So $F_0^{p,q}$ is a subsheaf of $F^{p,q}$ of full rank.
The compatibility of $\tau^0_{p,q}$ and $\tau_{p,q}$ implies
that $\sN_0^{p,q}$ is a subsheaf of
$$
\sN^{p,q}= {\rm Ker}(\tau_{p,q}: F^{p,q} \>>> F^{p-1,q+1}\otimes
\Omega^1_Y(\log S)).
$$
of maximal rank. Hence the induced morphism
$$
(\sN^{p,q})^\vee \to (\sN_0^{p,q})^\vee
$$
is an isomorphism over some dense open subset. By \ref{3.4},
iv), the sheaf $\sA\otimes F^{p,q}$ is a subsheaf of $E^{p,q}$
and by \ref{3.4}, i), the restriction $\theta_{p,q}|_{F^{p,q}}$
coincides with ${\rm id}_{\sA}\otimes \tau_{p,q}$. Using
the notations from \ref{3.4}, v), one obtains
\begin{equation}\label{3.7.1}
\sA\otimes \sN^{p,q}=\sA\otimes F^{p,q}\cap \sK^{p,q} \> \subset
>> \sK^{p,q}.
\end{equation}
By \ref{3.4}, v), the dual sheaf $(\sK^{p,q})^\vee$ is weakly
positive. (\ref{3.7.1}) induces morphisms
$$
(\sK^{p,q})^\vee \>>> \sA^{-1}\otimes (\sN^{p,q})^\vee
\>>> \sA^{-1}\otimes (\sN_0^{p,q})^\vee,
$$
surjective over some dense open subset, and we obtain \ref{3.7}.
\end{proof}

\begin{claim}\label{3.8} \ \ \
\begin{enumerate}
\item[i)]
If ${\rm Var}(f)=\dim(Y)$, then $(\sN_0^{p,q})^\vee$ is big.
\item[ii)] In general, for some $\alpha >0$ and for some
invertible sheaf $\lambda$ of Kodaira dimension $\kappa(\lambda)
\geq {\rm Var}(f)$ the sheaf
$$
S^\alpha((\sN_0^{p,q})^\vee)\otimes\lambda^{-1}
$$
is generically generated.
\end{enumerate}
\end{claim}
\begin{proof}
Let us consider as in \ref{2.9a} and \ref{3.1} some finite
morphism $\psi:Y'\to Y$, a desingularization $X'$ of $X\times_YY'$,
the induced morphisms $\psi':X'\to X$ and $f':X'\to Y'$,
$S'=\psi^{-1}(S)$, and $\Delta'={\psi'}^{-1}\Delta$.
In \ref{3.1} we constructed an injection
$$
{\psi}^*(F_0^{p,q}) \>\zeta_{p,q}>> {F'}_0^{p,q}= R^qf'_*(\Omega^p_{X'/Y'}(\log
\Delta') \otimes {\sL'}^{-1})/_{\rm torsion},
$$
compatible with the edge morphisms $\tau^0_{p,q}$ and ${\tau'}_{p,q}^0$.
Thereby we obtain an injection
\begin{equation*}
\psi^* \sN_0^{p,q}\>\zeta'_{p,q}>>
{\sN'}_0^{p,q}:= {\rm Ker}({\tau'}^0_{p,q}).
\end{equation*}
In case $X\times_YY'$ is non-singular, $\zeta_{p,q}$ and hence
$\zeta'_{p,q}$ are isomorphisms.

If ${\rm Var}(f)=\dim(Y)$, the conditions (\ref{3.5.1}) and
(\ref{3.5.2}) hold true. As in the proof of \ref{0.3}, iv),
there exists a finite covering $\psi: Y'\to Y$ with
$X\times_YY'$ non-singular, such that (\ref{3.3.1}) holds true
for the pullback family $X'\to Y'$. Since a sheaf is ample with
respect to some open set, if and only if it has the property on
some finite covering, we obtain the bigness of
$(\sN_0^{p,q})^\vee$ by applying \ref{3.7} to
$\psi^* (\sN_0^{p,q})^\vee=({\sN'}_0^{p,q})^\vee$.

In general \ref{kappa} allows to assume that $X \to Y$
fits into the diagram (\ref{1.7.1}) constructed in
\ref{1.7}. Let us write ${F_0^\#}^{p,q}$, and ${\sN^\#_0}^{p,q}$ for
the sheaves corresponding to $F_0^{p,q}$ and ${\sN^\#_0}^{p,q}$
on $Y^\#$ instead of $Y$.

As we have seen in \ref{3.2} the smoothness of $\eta$ implies
that $\eta^*{F^\#}^{p,q}={F'}^{p,q}$, and
\begin{equation}\label{3.8.1}
\eta^*{\sN^\#_0}^{p,q}\simeq {\sN'_0}^{p,q}\supset
\psi^*{\sN_0}^{p,q}.
\end{equation}
On $Y^\#$ we are in the situation where the variation
is maximal, hence i) holds true and the dual of the kernel
${\sN^\#}_0^{p,q}$ is big. So for any ample invertible sheaf $\sH$
we find some $\alpha>0$ and a morphism
$$
\bigoplus^r\sH\>>> S^\alpha(({\sN^\#}_0^{p,q})^\vee)
$$
which is surjective over some open set. Obviously the same holds
true for any invertible sheaf $\sH$, independent of the
ampleness. In particular we may choose for
any $\nu>1$ with $f_*\omega^\nu_{X/Y}\neq 0$ and for the number
$N_\nu$ given by \ref{1.7}, g) the sheaf
$$
\sH=\det(g^\#_*\omega^\nu_{Z^\#/Y^\#})^{N_\nu}.
$$
By \ref{1.7}, f) and by (\ref{3.8.1}), applied to $Y'\to
Y^\#$, the sheaf
\begin{multline*}
\eta^*(S^\alpha(({\sN^\#}_0^{p,q})^\vee)\otimes
\det(g^\#_*\omega^\nu_{Z^\#/Y^\#})^{-N_\nu})
= S^{\alpha}(({\sN'}_0^{p,q})^\vee)\otimes
\det(g_*\omega^\nu_{Z/Y'})^{-N_\nu} \\
\subset \psi^*(S^\alpha((\sN_0^{p,q})^\vee)\otimes \lambda_\nu^{-1})
\end{multline*}
is generically generated. By \ref{2.2} the same holds true for some power of
$$
S^\alpha((\sN_0^{p,q})^\vee)\otimes \lambda_\nu^{-1}.
$$
By \ref{2.5} and by the choice of $\lambda_\nu$ in \ref{1.7}, g), one finds
$\kappa(\lambda_\nu)\geq {\rm Var}(f)$ (see \ref{kappa}).
\end{proof}

To finish the proof of \ref{0.3}, i) and ii) we just have to
repeat the arguments used to prove \ref{3.5}, using \ref{3.6}.
By \ref{3.8} $\sO_Y=F_0^{n,0}$ can not lie in the kernel of
$\tau_{n,0}$. We choose $1 \leq m \leq n$ to be the largest
number with $\tau^m (F_0^{n,0})\neq 0.$
Then $\tau^m (F_0^{n,0})$ is contained in
$\sN_0^{n-m,m}\otimes S^m(\Omega^1_Y(\log S)),$
and we obtain morphisms of sheaves
\begin{equation}\label{3.8.2}
(\sN_0^{n-m,m})^\vee \>>> F_0^{n,0} \otimes
(\sN_0^{n-m,m})^\vee \> \neq 0 >> S^m(\Omega^1_Y(\log S)).
\end{equation}

Under the assumptions made in \ref{0.3}, ii) we take $\sP$ to be
the image of this morphism. By \ref{3.8}, i), this is the
image of a big sheaf, hence big.

If ${\rm Var}(f) < \dim(Y)$ \ref{3.8}, ii) implies that $S^\alpha(
(\sN_0^{n-m,m})^\vee)\otimes\lambda^{-1}$ is globally generated,
for some $\alpha >0$, and by (\ref{3.8.2}) one obtains a
non-trivial morphism
$$
\bigoplus^r \lambda \>>> S^{\alpha\cdot m}(\Omega_Y^1(\log S)).
$$
\qed

\section{Base spaces of families of smooth minimal models}
\label{4}

As promised in section one, we will show that in problem
\ref{0.2}, the bigness in b) follows from the weak positivity in
a). The corresponding result holds true for base
spaces of morphisms of maximal variation whose fibres are smooth
minimal models.

Throughout this section $Y$ denotes a projective manifold, and $S$ a
reduced normal crossing divisor in $Y$.

\begin{corollary}\label{4.1} Let $f:V\to U=Y\setminus S$
be a smooth family of $n$-dimensional projective manifolds with
${\rm Var}(f)=\dim(Y)$ and with $\omega_{V/U}$
$f$-semi-ample. If $\Omega^1_Y(\log S)$ is weakly positive, then
$\omega_Y(S)$ is big.
\end{corollary}

\begin{proof}
By \ref{0.3}, iii), there exists some $m>0$, a big coherent
subsheaf $\sP$, and an injective map
$$
\sP \> \subset >> S^m(\Omega^1_Y(\log S)).
$$
Its cokernel $\sC$, as the quotient of a weakly positive sheaf, is
weakly positive, hence $\det(S^m(\Omega^1_Y(\log S)))$ is the
tensor product of the big sheaf $\det(\sP)$ with the weakly
positive sheaf $\det(\sC)$.
\end{proof}
\begin{corollary}[Kov\'acs, \cite{Kov3}, for
$S=\emptyset$]\label{4.2}\ \\
If $T_Y(-\log S)$ is weakly positive, then there exists
\begin{enumerate}
\item[a)] no non-isotrivial smooth projective family $f:V\to U$
of canonically polarized manifolds.
\item[b)] no smooth projective family $f:V \to U$ with
${\rm Var}(f)=\dim(U)$ and with $\omega_{V/U}$ $f$-semi-ample.
\end{enumerate}
\end{corollary}
\begin{proof}
In both cases \ref{0.3} would imply for some $m>0$ that
$S^m(\Omega^1_Y(\log S))$ has a subsheaf $\sA$ of positive
Kodaira dimension. But $\sA^\vee$, as a quotient of a weakly
positive sheaf, must be weakly positive, contradicting $\kappa(\sA)>0$.
\end{proof}
There are other examples of varieties $U$ for which
$S^m(\Omega^1_Y(\log S))$ can not contain a subsheaf of strictly
positive Kodaira dimension or more general, for which
\begin{equation}\label{4.2.1}
H^0(Y,S^m(\Omega^1_Y(\log S)))=0 \mbox{ \ \ \ for all \ \ \ }
m>0.
\end{equation}
The argument used in \ref{4.2} carries over and excludes the
existence of families, as in \ref{4.2}, a) or b). For example,
(\ref{4.2.1}) has been verified by Br\"uckmann for $U=H$ a complete
intersection in $\P^N$ of codimension $\ell < \frac{N}{2}$ (see
\cite{B-R} for example). As a second application of this result,
one can exclude certain discriminant loci for families of canonically
polarized manifolds in $\P^N$. If $H=H_1+\cdots +H_\ell$ is a
normal crossing divisor in $\P^N$, and $\ell < \frac{N}{2}$,
then for $U=\P^N\setminus H$ the conclusions in
\ref{4.2} hold true. In order allow a proof by induction, we
formulate both results in a slightly more general setup.
\begin{corollary}\label{4.3}
For $\ell < \frac{N}{2}$ let $H=H_1+\cdots +H_\ell$ be a normal
crossing divisor in $\P^N$. For $0\leq r\leq l$ define
\begin{gather*}
H = \bigcap_{j=r+1}^\ell H_j, \ \ \ S_i = H_i|_H, \ \ \
S = \sum_{i=1}^r S_i,
\mbox{ \ \ and \ \ } U= H \setminus S
\end{gather*}
(where for $l=r$ the intersection with empty index set is
$H=\P^N$). Then there exists
\begin{enumerate}
\item[a)] no non-isotrivial smooth projective family $f:V\to U$
of canonically polarized manifolds.
\item[b)] no smooth projective family $f:V \to U$ with
${\rm Var}(f)=\dim(U)$ and with $\omega_{V/U}$ $f$-semi-ample.
\end{enumerate}
\end{corollary}
\begin{proof}
Let $\sA$ be an invertible sheaf of Kodaira dimension
$\kappa(\sA) > 0$. Replacing $\sA$ by some power, we may assume
that $\dim(H^0(H,\sA)) > r+1 $. We have to verify that
there is no injection $\sA \to S^m(\Omega^1_{H}(\log S)).$

For $r=0$ such an injection would contradict the vanishing
(\ref{4.2.1}) shown in \cite{B-R}. Hence starting with $r=0$
we will show the non-existence of the subsheaf $\sA$
by induction on $\dim(H)=N-\ell+r$ and on $r$.

The exact sequence
$$
0 \to \Omega^1_{S_{r}}(\log (S_{1} + \cdots +
S_{r-1})) \to
\Omega^1_{H}(\log S)|_{S_{r}} \to
\sO_{S_{r}} \to 0
$$
induces a filtration on
$S^m(\Omega^1_{H}(\log S))|_{S_{r}}$
with subsequent quotients
$$
S^\mu(\Omega^1_{S_{r}}(\log (S_{1} + \cdots +
S_{r-1})))
$$
for $\mu = 0, \cdots , m$. By induction none of those quotients
can contain an invertible subsheaf of positive Kodaira
dimension. Hence either the restriction of $\sA$ to $S_r$ is a
sheaf with $\kappa(\sA|_{S_r}) \leq 0$, hence
$\dim(H,\sA(-S_r))> r$, or the image of $\sA$ in
$S^m(\Omega^1_{H}(\log S))|_{S_{r}}$
is zero. In both cases
$$
S^m(\Omega^1_{H}(\log (S))\otimes
\sO_H(-{S_{r}})
$$
contains an invertible subsheaf $\sA_1$ with at least two
linearly independent sections, hence of positive Kodaira
dimension. Now we repeat the same argument a second time:
$$
(S^m(\Omega^1_{H}(\log S))\otimes
\sO_H(-{S_{r}}))|_{S_{r}}
$$
has a filtration with subsequent quotients
$$
(S^\mu(\Omega^1_{S_{r}}(\log (S_{1} + \cdots +
S_{r-1}))\otimes \sO_H(-{S_{r}}))|_{S_{r}}).
$$
$\sO_H(S_r)$ is ample, hence by induction
none of those quotients can have a non-trivial section. Repeating this
argument $m$ times, we find an invertible sheaf contained in
\begin{multline*}
S^m(\Omega^1_{H}(\log S))\otimes
\sO_H(-m\cdot {S_{r}})= S^m(\Omega^1_{H}(\log S)\otimes
\sO_H(-{S_{r}})) \\
\subset S^m(\Omega^1_{H}(\log (S_{1} + \cdots + S_{r-1}))),
\end{multline*}
and of positive Kodaira dimension, contradicting the induction hypothesis.
\end{proof}

\section{Subschemes of moduli stacks of canonically polarized manifolds}
\label{5}

Let $M_h$ denote the moduli scheme of canonically polarized
$n$-dimensional manifolds with Hilbert polynomial $h$.
In this section we want to apply \ref{0.3} to obtain properties
of submanifolds of the moduli stack. Most of those remain
true for base spaces of smooth families with a relatively semi-ample
dualizing sheaf, and of maximal variation.
\begin{assumptions}\label{5.1} Let $Y$ be a projective manifold,
$S$ a normal crossing divisor and $U=Y-S$. Consider the
following three setups:
\begin{enumerate}
\item[a)] There exists a quasi-finite morphism $\varphi:U \to M_h$
which is induced by a smooth family $f:V\to U$ of canonically
polarized manifolds.
\item[b)] There exists a smooth family $f:V\to U$ with
$\omega_{V/U}$ $f$-semi-ample and with ${\rm Var}(f)=\dim(U)$.
\item[c)] There exists a smooth family $f:V\to U$ with
$\omega_{V/U}$ $f$-semi-ample and some $\nu \geq 2$ for which
the following holds true. Given a non-singular
projective manifold $Y'$, a normal crossing divisors $S'$ in $Y'$,
and a quasi-finite morphism $\psi':U'=Y'\setminus S' \>>> U$, let
$X'$ be a non-singular projective compactification of
$V\times_UU'$ such that the second projection induces a morphism
$f':X'\to Y'$. Then the sheaf $\det(f'_*\omega_{X'/Y'}^\nu)$ is
ample with respect to $U'$.
\end{enumerate}
\end{assumptions}
Although we are mainly interested in the cases \ref{5.1}, a) and
b), we included the quite technical condition c), since this is
what we really need in the proofs.

The assumption made in a) implies the one in c). In fact, if
$\varphi$ is quasi-finite, the same holds true for
$\varphi\circ\psi':U'\to M_h$, and by \ref{2.6}, iii),
$\det(f'_*\omega_{X'/Y'}^\nu)$ is ample with respect to $U'$.

Under the assumption b), it might happen that we have to
replace $U$ in c) by some smaller open subset $\tilde{U}$.
To this aim start with the open set $U_g$ considered in
\ref{1.6}. Applying \ref{2.6}, i), to the family $g:Z \to Y'$
in (\ref{1.5.1}), one finds an open subset $\tilde{U}'$
of $Y'$ with $g_*\omega^\nu_{Z/Y'}$ ample with respect to
$\tilde{U}'$. We may assume, of course, that $\tilde{U}'$
is the preimage of $\tilde{U}\subset U_g$.
Since $g:Z \to Y'$ is birational to a mild morphism over $Y'$,
the same holds true for all larger coverings, and the condition c)
follows by flat base change, for $\tilde{U}$ instead of $U$.

Let us start with a finiteness result for morphisms from
curves to $M_h$, close in spirit to the one obtained in
\cite{B-V}, 4.3, in case that $M_h$ is the moduli space of
surfaces of general type. Let $C$ be a projective non-singular
curve and let $C_0$ be a dense open subset of $C$. By \cite{Gro} the morphisms
$$
\pi : C \to Y \mbox{ \ \  with \ \ } \pi (C_0) \subset U
$$
are parameterized by a scheme ${\rm \bf H}:={\rm \bf Hom}
((C,C_0),(Y,U))$, locally of finite type.
\begin{theorem}\label{5.2}\ \ \
\begin{enumerate}
\item[i)] Under the assumptions made in \ref{5.1}, a) or c), the
scheme ${\rm \bf H}$ is of finite type.
\item[ii)] Under the assumption \ref{5.1}, b), there exists an
open subscheme $U_g$ in $U$ such that there are only finitely many
irreducible components of ${\rm \bf H}$ which parameterize morphisms
$\pi : C \to Y$ with $\pi (C_0) \subset U$ and $\pi (C_0) \cap
U_g \neq \emptyset.$
\end{enumerate}
\end{theorem}
\begin{proof}
Let us return to the notations
introduced in \ref{1.5} and \ref{1.6}. There we considered
an open dense non-singular subvariety $U_g$ of $U$, depending on the
construction of the diagram (\ref{1.5.1}). In particular $U_g$
embeds to $Y$.

Let $\sH$ be an ample invertible sheaf on $Y$.
In order to prove i) we have to find some constant $c$
which is an upper bound for $\deg(\pi^*\sH)$, for all
morphisms $\pi: C \to Y$ with $\pi(C_0)\subset U$.
For part ii) we have to show the same, under the additional
assumption that $\pi(C)\cup U_g \neq \emptyset$.
Let us start with the latter.

By \ref{2.6}, iii) in case \ref{5.1}, a), or by
assumption in case \ref{5.1}, c), one finds the sheaf $\lambda_\nu$, defined
in \ref{1.5}, d), to be ample with respect to $U_g$. For part
ii) of \ref{5.2}, i.e. if one just assumes that ${\rm
Var}(f)=\dim(U)$, we may use \ref{2.5}, and choose $U_g$ a bit
smaller to guarantee the ampleness of $g_*\omega_{Z/Y'}^\nu$
over $\psi^{-1}(U_g)$.

Replacing $N_\nu$ by some multiple and $\lambda_\nu$ by some
tensor power, we may assume that $\lambda_\nu\otimes\sH^{-1}$ is
generated by global sections over $U_g$.

Assume first that $\pi(C_0)\cap U_g \neq \emptyset$.
Let $h: W\to C$ be a morphism between projective manifolds,
obtained as a compactification of $X\times_YC_0 \to C_0$. By
definition $h$ is smooth over $C_0$. In \ref{1.6} we have
shown, that
$$
\deg(\pi^*\lambda_\nu) \leq N_\nu\cdot
\deg(\det(h_*\omega_{W/C}^\nu)).
$$
On the other hand, upper bounds for the right hand side have
been obtained for case a) in \cite{B-V}, \cite{Kov2} and in
general in \cite{V-Z}. Using the notations from \cite{V-Z},
$$
\deg(\det(h_*\omega_{W/C}^\nu)) \leq
(n\cdot(2g(C)-2+s)+s)\cdot\nu\cdot {\rm
rank}(h_*\omega_{W/C}^\nu)\cdot e,
$$
where $g(C)$ is the genus of $C$, where $s=\#(C-C_0)$, and where
$e$ is a positive constant, depending on the general fibre of $h$.
In fact, if $F$ is a general fibre of $h$, the constant $e$ can be
chosen to be $e(\omega_F^\nu)$. Since the latter is upper
semicontinous in smooth families (see \cite{E-V} or \cite{Vie},
5.17) there exists some $e$ which works for all possible curves.
Altogether, we found an upper bound for $\deg(\pi^*\sH)$,
whenever the image $\pi(C)$ meets the dense open subset $U_g$ of $U$.

In i), the assumptions made in \ref{5.1}, a) and c) are
compatible with restriction to subvarieties of $U$, and we may
assume by induction, that we already obtained similar bounds for
all curves $C$ with $\pi(C_0)\subset (U\setminus U_g)$.
\end{proof}

From now on we fix again a projective non-singular
compactification with $S=Y\setminus U$ a normal crossing
divisor. Even if $\varphi:U\to M_h$ is quasi finite, one can not
expect $\Omega^1_Y(\log S)$ to be ample with respect to $U$,
except for $n=1$, i.e. for moduli of curves.
For $n >1$ there are obvious counter examples.
\begin{example}\label{5.3}
Let $ g_1: Z_1\to C_1 $ and $ g_2:Z_2\to C_2$ be two
non-isotrivial families of curves over curves $ C_1$ and $ C_2,$
with degeneration loci $S_1$ and $S_2$, respectively. We assume
both families to be semi-stable, of different genus, and we
consider the product
$$
f: X=Z_1\times Z_2\>>> C_1\times C_2=Y,
$$
the projections $ p_i : Y\to C_i$, and
the discriminant locus $S=p^{-1}_1(S_1)\cup p^{-1}_2(S_2).$
For two invertible sheaves $\sL_i$ on $C_i$ we write
$$
\sL_1\boxplus\sL_2=p_1^*\sL_1\oplus p_2^*\sL_2 \mbox{ \ \ and \
\ } \sL_1\boxtimes\sL_2=p_1^*\sL_1\otimes p_2^*\sL_2.
$$
For example
$$
S^2(\sL_1 \boxplus \sL_2) = p_1^* \sL_1^2 \oplus p_2^*\sL_2^2 \oplus
\sL_1\boxtimes \sL_2.
$$

The family $f$ is non-isotrivial, and it induces a
generically finite morphism to the moduli space of surfaces of
general type $M_h$, for some $h$. Obviously,
$$
\Omega^1_Y(\log S)= \Omega^1_{C_1}(\log S_1)\boxplus
\Omega^1_{C_2}(\log S_2):= p_1^*(\Omega^1_{C_1}(\log S_1))\oplus
p_2^*(\Omega^1_{C_2}(\log S_2))
$$
can not be ample with respect to any open dense subset.

Let us look, how the edge morphisms $\tau_{p,q}$ defined in
section \ref{3} look like in this special case. To avoid conflicting
notations, we write $G_i^{p,q}$ instead of $F^{p,q}$, for the
two families of curves, and
$$
\sigma_i:{G_i}^{1,0} = {g_i}_* \sO_{Z_i}=\sO_{C_i} \>>> {G_i}^{0,1}
\otimes \Omega^1_{C_i}(\log S_i)
$$
for the edge morphisms. The morphism
$$
\tau^2=\tau_{1,1}\circ\tau_{2,0}: F^{2,0}=\sO_Y \>>>
F^{0,2}\otimes S^2(\Omega_Y^1(\log S)),
$$
considered in the proof of \ref{0.3}, i) and ii), thereby
induces three maps,
$$
t_i:F^{0,2\vee}\to S^2(p_i^*\Omega^1_{C_i}(\log S_i)),
$$
for $i=1, \ 2$, and
$$
t: F^{0,2\vee}\>>> \Omega^1_{C_1}(\log S_1)\boxtimes
\Omega^1_{C_2}(\log S_2).
$$
Since $F^{0,2\vee}=g_{1*}\omega^2_{Z_1/C_1}\boxtimes
g_{2*}\omega^2_{Z_2/C_2}$ is ample the first two morphisms
$t_1$ and $t_2$ must be zero.
$$
F^{1,1} = R^1f_*((\omega_{Z_1/C_1}\boxplus
\omega_{Z_2/C_2}) \otimes \omega_{X/Y}^{-1})
\simeq R^1f_* (\omega_{Z_1/C_1}^{-1}\boxplus \omega_{Z_2/C_2}^{-1}),
$$
where the isomorphism interchanges the two factors. In particular
$$
F^{1,1} = G^{0,1}_1\boxplus G^{0,1}_2,
$$
and one has $\tau_{2,0} = \sigma_1 \boxplus \sigma_2$.
Its image lies in the direct factor
$$
G':=G^{0,1}_1\otimes \Omega^1_{C_1}(\log S_1) \boxplus G^{0,1}_2 \otimes
\Omega^1_{C_2}(\log S_2)
$$
of $F^{1,1}\otimes \Omega^1_Y(\log S)$. The picture should
be the following one:
$$
F^{0,2} \otimes \Omega^1_Y(\log S) =
(G_1^{0,1}\boxtimes G_2^{0,1})\otimes (\Omega^1_{C_1}(\log
S_1)\boxplus \Omega^1_{C_2}(\log S_2)),
$$
and $\tau_{1,1}|_{G^{0,1}_1}={\rm id}_{G^{0,1}_1}\otimes
\sigma_2$ with image in
$$
(G_1^{0,1}\boxtimes G_2^{0,1})\otimes \Omega^1_{C_2}(\log S_2).
$$
Hence $\tau_{1,1}\circ \tau_{2,0}$ is the sum of the two maps
$\tau_{1,1}|_{G^{0,1}_i} \circ \sigma_i$, both with image in
$$
(G_1^{0,1}\boxtimes G_2^{0,1})\boxtimes \Omega^1_{C_1}(\log
S_1)\otimes \Omega^1_{C_2}(\log S_2).
$$
\end{example}
In general, when there exists a generically finite morphism
$\varphi:C_1\times C_2\to M_h$ induced by $f:V\to U$, the
picture should be quite similar, however we were unable to
translate this back to properties of the general
fibre of $f$. However, for moduli of surfaces there can not
exist a generically finite morphism from the product of three
curves. More generally one obtains from \ref{0.3}:

\begin{corollary}\label{5.4} Let $U=C^0_1 \times \cdots
C^0_\ell$ be the product of $\ell$ quasi-projective curves, and
assume there exists a smooth family $f: V\to U$ with
$\omega_{V/U}$ $f$-semi-ample and with ${\rm Var}(f)=\dim(U)$.
Then $\ell \leq n=\dim(V)-\dim(U)$.
\end{corollary}
\begin{proof}
For $C_i$, the non-singular compactification of $C^0_i$, and for
$S_i=C_i\setminus C_i^0$, a compactification of $U$ is given by
$Y=C_1\times \cdots \times C_\ell$ with boundary divisor
$S=\sum_{i=1}^\ell pr_i^* S_i$. Then
$$
S^m(\Omega_Y^1(\log S))=\bigoplus
S^{j_1}(pr_1^*\Omega^1_{C_1}(\log S_1))\otimes \cdots \otimes
S^{j_\ell}(pr_\ell^*\Omega^1_{C_\ell}(\log S_\ell))
$$
where the sum is taken over all tuples $j_1, \ldots ,j_\ell$
with $j_1+\cdots+j_\ell=m$. If $\ell > m$, each of the factors
is the pullback of some sheaf on a strictly lower
dimensional product of curves, hence for $\ell > m$ any morphism
from a big sheaf $\sP$ to $S^m(\Omega_Y^1(\log S))$ must be trivial.
If $\psi:Y' \to Y$ is a finite covering, the same holds true
for $\psi^*S^m(\Omega_Y^1(\log S))$. By \ref{0.3}, iv), there exists
such a covering, some $m \leq n$ and a big subsheaf of
$\psi^*S^m(\Omega_Y^1(\log S))$, hence $\ell \leq m \leq n$.
\end{proof}

The next application of \ref{0.3} is the rigidity of generic
curves in moduli stacks. If $\varphi:U\to M_h$ is induced by a family,
\ref{0.3}, ii) provides us with a big subsheaf $\sP$ of
$S^m(\Omega_Y^1(\log S))$, and if we do not insist that $m\leq
n$, the same holds true whenever there exists a family $V\to U$,
as in \ref{5.1}, b). In both cases, replacing $m$ by some
multiple, we find an ample invertible sheaf $\sH$ on $Y$ and an
injection
$$
\iota: \bigoplus \sH \>>> S^m(\Omega_Y^1(\log S)).
$$
Let $U_1$ be an open dense subset in $U$, on which $\iota$
defines a subbundle.
\begin{corollary}\label{5.5}
Under the assumption a) or b) in \ref{5.1} there exists an open
dense subset $U_1$ and for each point $y\in U_1$ a curve
$C_0\subset U$, passing through $y$, which is rigid, i.e.:
If for a reduced curve $T_0$ and for $t\in T_0$
there exists a morphism $\rho: T_0 \times C_0 \to U$ with
$\rho({\{t_0\}\times C_0})={C_0}$, then $\rho$ factors through
$pr_2:T_0 \times C_0 \to C_0$.
\end{corollary}
\begin{proof}
Let $\pi:\P=\P(\Omega_Y^1(\log S))\to Y$ be the projective
bundle.
$$
\iota:\bigoplus \sH \>>> \pi_*\sO_\P(m)
$$
defines sections of $\sO_Y(m)\otimes \pi^*\sH^{-1}$, which
are not all identically zero on $\pi^{-1}(y)$ for $y\in U_1$.
Hence there exists a non-singular curve $C_0\subset U$ passing through
$y$, such that the composite
\begin{equation}\label{5.5.1}
\bigoplus \sH|_U \>>> S^m(\Omega_U^1) \>>>
S^m(\Omega_{C_0}^1)
\end{equation}
is surjective over a neighborhood of $y$. For a
nonsingular curve $T_0$ and $t_0\in T_0$ consider a morphism
$\phi_0:T_0\times C_0 \to U$, with $\phi_0({\{t_0\}\times C_0})=
{C_0}$. Let $T$ and $C$ be projective non-singular curves,
containing $T_0$ and $C_0$ as the complement of divisors $\Theta$
and $\Gamma$, respectively. On the complement $W$ of a
codimension two subset of $T\times C$ the morphism
$\phi_0$ extends to $\phi:W \to Y$.
Then $\iota$ induces a morphism
$$
\phi^*\bigoplus \sH \>>> \phi^* S^m(\Omega_Y^1(\log S)) \>>>
S^m(\Omega_T^1(\log \Theta)\boxplus \Omega_C^1(\log \Gamma))|_W
$$
whose composite with
\begin{multline*}
S^m(\Omega_T^1(\log \Theta)\boxplus \Omega_C^1(\log
\Gamma))|_W\>>> S^m(pr_2^* \Omega_C^1(\log \Gamma))|_W\\
\>>> S^m(\Omega_{\{t\}\times C}(\log \Gamma))|_{W\cap \{t\}\times C}
\end{multline*}
is non-zero for all $t$ in an open neighborhood of $t_0$,
hence
$$
\phi^*\bigoplus \sH \>>>
pr_2^* S^m(\Omega_C^1(\log \Gamma))|_{W}
$$
is surjective over some open dense subset. Since $\sH$ is ample,
this is only possible if $\phi:W \to Y$ factors
through the second projection $W \to T\times C \to C$.
\end{proof}
Assume we know in \ref{5.5} that $\Omega^1_Y(\log S)$ is ample
over some dense open subscheme $U_2$. Then the morphism
(\ref{5.5.1}) is non-trivial for all curves $C_0$ meeting $U_2$,
hence the argument used in the proof of \ref{5.5} implies, that a
morphisms $\pi: C_0 \to U$ with $\pi(C_0)\cap U_2\neq \emptyset$
has to be rigid.
If $U_2=U$, this, together with \ref{5.2} proves the next
corollary.
\begin{corollary}\label{5.6}
Assume in \ref{5.2}, i), that $\Omega_Y^1(\log S)$ is ample with
respect to $U$. Then $\rm \bf H$ is a finite set of points.
\end{corollary}
The generic rigidity in \ref{5.5}, together with the finiteness result in
\ref{5.2}, implies that subvarieties of the moduli stacks have a
finite group of automorphism. Again, a similar statement holds
true under the assumption \ref{5.1}, b).
\begin{theorem}\label{5.7}
Under the assumption \ref{5.1}, a) or b) the automorphism group
${\rm Aut}(U)$ of $U$ is finite.
\end{theorem}
\begin{proof}
Assume ${\rm Aut}(U)$ is infinite, and choose an infinite countable
subgroup
$$G \subset {\rm Aut}(U).$$
Let $U_1$ be the open subset
of $U$ considered in \ref{5.5}, and let $U_g$ be the open subset
from \ref{5.2}, b). We may assume that $U_1 \subset U_g$ and write
$\Gamma= U\setminus U_1$. Since
$$
\bigcup_{g\in G} g(\Gamma) \neq U
$$
we can find a point $y\in U_1$ whose $G$-orbit is an infinite
set contained in $U_1$. By \ref{5.5} there are rigid smooth curves
$C_0 \subset U$ passing through $y$. Obviously, for all $g\in
G$ the curve $g(C_0)\subset U$ is again rigid, it meets $U_1$,
hence $U_g$ and the set of those curves is infinite,
contradicting \ref{5.2}, b).
\end{proof}
\section{A vanishing theorem for sections of symmetric powers of
logarithmic one forms}\label{7}

\begin{proposition} \label{7.1}
Let $Y$ be a projective manifold and let $D = D_1 + \ldots +
D_r$ and $S = S_1 + \ldots + S_{\ell}$ be two reduced divisors
without common component. Assume that
\begin{enumerate}
\item[i)] $S + D$ is a normal crossing divisor.
\item[ii)] For no subset $J \subseteq \{ 1, \ldots , \ell\}$ the
intersection
$$
S_J = \bigcap_{j \in J} S_i
$$
is zero dimensional.
\item[iii)] $T_Y ( - {\rm log} \ D) = (\Omega^{1}_{Y} ({\rm log}
\ D))^{\vee}$ is weakly positive over $U_1 = Y - D$.
\end{enumerate}
Then for all ample invertible sheaves $\sA$ and for all $m \geq
1$
$$
H^0 (Y, S^m (\Omega^{1}_{Y} ({\rm log} (D + S))) \otimes
\sA^{-1} ) =0.
$$
\end{proposition}

\begin{corollary} \label{7.2}
Under the assumption i), ii) and iii) in \ref{7.1} there exists
no smooth family $f : V \to U = Y \setminus (S +D)$ with ${\rm
Var} (f) = \dim Y$. In particular there is no generically finite
morphism $U \to M_h$, induced by a family.
\end{corollary}

\begin{proof}
By \ref{0.3}, iii) the existence of such a family implies that
for some $m >0$ the sheaf $S^m (\Omega^{1}_{Y} ({\rm log}
(D+S)))$ contains a big coherent subsheaf $\sP$. Replacing $m$
by some multiple, one can assume that $\sP$ is ample and
invertible, contradicting \ref{7.1}.
\end{proof}

\noindent
{\it Proof of \ref{cor5} and of the second part of \ref{cor3}.}
\ Since for an abelian variety $Y$ the sheaf $\Omega^{1}_{Y}$ is
trivial, and since the condition ii) in \ref{7.1} is obvious for
$\ell < \dim Y$, part a) of \ref{cor5} is a special case of
\ref{7.2}.

For b) again i) and ii) hold true by assumption. For iii) we
remark, that
$$
\Omega^{1}_{\P^{\nu_i}} ({\rm log} \ D^{(\nu_i)} ) =
\oplus^{\nu_i} \sO_{\P^{\nu_i}},
$$
hence $\Omega^{1}_{Y} ({\rm log} \ D)$ is again a direct sum of
copies of $\sO_Y$. Assume that $\sA$ is an invertible subsheaf of
$S^m (\Omega^{1}_{Y} ({\rm log} (D+S)))$, for some $m > 0$. If
$\kappa (\sA) > 0$, then for some $\mu_i \in \N$,
$$
\sA = \sO_Y (\mu_1, \ldots , \mu_k ) = \bigotimes^{k}_{i=1}
pr^{*}_{i} \sO_{\P^{\nu_i}} (\mu_i).
$$
By \ref{7.1} not all the $\mu_i$ can be strictly larger than
zero, hence
$$
\kappa (\sA) \leq {\rm Max} \{ \dim (Y) - \nu_{i} ; \
i=1, \ldots , k\} = M.
$$
By \ref{0.3}, i), for any morphism $\varphi : U \to M_{h}$,
induced by a family $f : V \to U$, one has
$$
{\rm Var} (f) = \dim (\varphi (U)) \leq M.
$$
\ref{0.3}, iii), implies that there exists no smooth family $f:
V \to U$ of maximal variation, and with $\omega_{V/U}$
$f$-semi-ample. \qed

\begin{remark} \label{7.3}
In \ref{cor5}, a), one can also show, that for an abelian
variety $Y$ and for a morphism $\varphi : Y \to M_h$, induced by
a family,
$$
\dim (\varphi (Y)) \leq \dim (Y) - \nu,
$$
where $\nu$ is the dimension of the smallest simple abelian
subvariety of $Y$. In fact, by the Poincar\'e decomposition
theorem, $Y$ is isogenous to the product of simple abelian
varieties, hence replacing $Y$ by an \'etale covering, we may
assume that
$$
Y = Y_1 \times \ldots \times Y_k
$$
with $Y_i$ simple abelian, and with
$$
\nu = \dim (Y_1) \leq \dim (Y_2) \leq \ldots \leq \dim Y_k.
$$
Since an invertible sheaf of positive Kodaira dimension on a
simple abelian variety must be ample, one finds that for an
non-ample invertible sheaf $\sA$ on $Y$
$$
\kappa (\sA) \leq \dim (Y) - \nu.
$$
\end{remark}

Before proving \ref{7.1} let us show that for $Y = \P^2$ we can
not allow $S$ to have two irreducible components of high degree,
even if $D = 0$.

\begin{example} \label{7.4}
Given a surface $Y$ and a normal crossing divisor $S+D$, with $S
= \sum^{\ell}_{i=1} S_i$, consider the two exact sequences
$$
0 \to \sO_Y (-S_i) \to \sO_Y \to \sO_{S_i} \to 0 \ \mbox{ \ \ \
\ \ \ and}
$$
$$
0 \to \Omega^{1}_{Y} ({\rm log} \ D) \to \Omega^{1}_{Y} ({\rm
log} (D+S)) \to \bigoplus^{\ell}_{i=1} \sO_{S_i} \to 0.
$$
Writing $c (\sE) =1 + c_1 (\sE) + c_2 (\sE)$ for a sheaf $\sE$
on $Y$, one finds
$$
c (\sO_{S_i}) = 1 + S_i + S^{2}_{i}
$$
and
$$
c (\Omega^{1}_{Y} ({\rm log} (D+S))) = c (\Omega^{1}_{Y} ({\rm
log} \ D)) \cdot \prod^{\ell}_{i=1} (1 + S_i + S^{2}_{i} ).
$$
Hence
$$
c_2 (\Omega^{1}_{Y} ({\rm log} (D + S))) = c_2 (\Omega^{1}_{Y}
({\rm log} \ D)) + \sum^{\ell}_{i=1} c_1 (\Omega^{1}_{Y} ({\rm
log} \ D)). S_i + \sum_{i < j} S_i . S_j + \sum^{\ell}_{i=1}
S^{2}_{i}
$$
and $c_1 (\Omega^{1}_{Y} ({\rm log} (D+S))) = c_1
(\Omega^{1}_{Y}) + D+S$.
The Riemann-Roch theorem for vector bundles on surfaces and the
isomorphism
$$
S^m (\Omega^{1}_{Y} ({\rm log} (D+S))^{\vee}) \otimes \omega_Y =
S^m (\Omega^{1}_{Y} ({\rm log} (D+S))) \otimes \omega_{Y}
\otimes \omega_Y (D+S)^{-m}
$$
imply that for an invertible sheaf $\sA$
\begin{multline*}
h^0 (Y, S^m (\Omega^{1}_{Y} ({\rm log} (D+S))) \otimes \sA^{-1}
)\\ + h^0 (Y, S^m (\Omega^{1}_{Y} ({\rm log} (D+S))) \otimes
\omega_Y \otimes \omega_Y (D+S)^{-m} \otimes \sA)\\
\geq \frac{m^3}{6} (c_1 (\Omega^{1}_{Y} ({\rm log} (D+S)))^2 -
c_2 (\Omega^{1}_{Y} ({\rm log} (D+S)))) + O(m^2),
\end{multline*}
where $O(m^2)$ is a sum of terms of order $\leq 2$ in $m$.
If $\omega_Y (D+S)$ is big and if
\begin{equation} \label{7.4.1}
c_1 (\Omega^{1}_{Y} ({\rm log} (D+S)))^2 > c_2 (\Omega^{1}_{Y}
({\rm log} (D+S)))
\end{equation}
then for $m \gg 0$ the sheaf $\omega_Y \otimes \omega_Y
(D+S)^{-m} \otimes \sA$ is a subsheaf of $\sA^{-1}$ and
$$
h^0 (Y, S^m (\Omega^{1}_{Y} ({\rm log} (D+S))) \otimes \sA^{-1}
) \neq 0.
$$
For $Y = \P^2$ and for a coordinate system $D = D_0 + D_1 + D_2$,
\begin{multline*}
c_1 (\Omega^{1}_{Y} ({\rm log} (D+S)))^2 - c_2 (\Omega^{1}_{Y}
({\rm log} \ D +S)))\\ = ( \sum^{\ell}_{i=1} S_i )^2 -
\sum^{\ell}_{i=1} S^{2}_{i} - \sum_{i<j} S_i . S_j = \sum_{i
< j} S_i . S_j,
\end{multline*}
and as soon as $S$ has more than one component, (\ref{7.4.1})
holds true. So in \ref{7.1}, for $Y=\P^2$ and $\ell >1$, the arguments
used to proof \ref{cor5} fail.

We do not know, whether $U = \P^2 \setminus (S_1 +
S_2)$ can be the base of a non-isotrivial family of canonically
polarized manifolds.

For $D = 0$ and $Y = \P^2$, one finds
\begin{multline*}
c_1 (\Omega^{1}_{Y} ({\rm log} \ S))^2 - c_2 (\Omega^{1}_{Y}
({\rm log} \ S))= 6 - 3 \cdot {\rm deg} (S) + \sum_{i<j} S_i
. S_j\\
= 3\cdot(2-\deg(S)) + \sum_{i<j} \deg(S_i)\cdot \deg(S_j).
\end{multline*}
Assume that $2 \leq {\rm deg} (S_1) \leq {\rm deg} (S_2) \ldots
\leq {\rm deg} (S_{\ell})$. Then the only cases where
(\ref{7.4.1}) does not hold true are $\ell = 2$ and ${\rm deg} (S_1) =
2$, or $\ell =3$ and ${\rm deg} (S_i) =2$ for $i=1, 2, 3$.
Again we do not know any example of a non-isotrivial family
over $U=\P^2\setminus S$.
\end{example}

As a first step in the proof of \ref{7.1} we need
\begin{lemma} \label{7.5}
Let $\sE$ and $\sF$ be locally free sheaves on $Y$. Assume that,
for a non-singular divisor $B$, for some ample invertible sheaf
$\sA$, and for all $m \geq 0$
$$
H^0 (Y, S^m (\sF) \otimes \sA^{-1}) = H^0 (B, S^m (\sF) \otimes
\sA^{-1} |_B ) =0.
$$
Assume moreover that there exists an exact sequence
$$
0 \to \sF \to \sE \to \sO_B \to 0.
$$
Then for all $m > 0$
$$
H^0 (Y, S^m (\sE) \otimes \sA^{-1}) =0.
$$
\end{lemma}

\begin{proof}
Write $\pi : \P = \P (\sE) \to Y$. The surjection $\sE \to
\sO_B$ defines a morphism $s : B \to \P$. For the ideal $I$ of
$s(B)$ the induced morphism $\pi^* \sF \to I \otimes \sO_{\P}
(1)$ is surjective, as well as the composite
$$
\tilde{\pi}^* \sF \to \delta^* (I \otimes \sO_{\P} (1)) \to
\sO_{\tilde{\P}} (-E) \otimes \delta^* \sO_{\P} (1),
$$
where $\delta : \tilde{\P} \to \P$ is the blowing up of $I$ with
exceptional divisor $E$, and where $\tilde{\pi} = \pi \circ
\delta$.

Let us write $M+1$ for the rank of $\sE$. For $y \in B$ and
$p=s(y)$ let
$$
\delta_y:\tilde{\P}_y \to \P^M=\pi^{-1}(y)
$$
be the blowing up of $p$, with exceptional divisor $F$.
Then
$$
\tilde{\pi}^{-1}(y) = \tilde{\P}_y \cup \P^M \mbox{ \ \ with \ \
} F=\tilde{\P}_y \cap \P^M.
$$
In particular $\tilde{\pi}$ is equidimensional, hence flat.
For $0\leq \mu \leq m$ and for $i>0$
$$
H^i(\tilde{\P}_y, \sO_{\tilde{\P}_y}(-(\mu+1)\cdot
F)\otimes\delta_y^* \sO_{\P^M}(m))=0
$$
and
$$
H^0(\tilde{\P}_y, \sO_{\tilde{\P}_y}(-(\mu+1)\cdot
F)\otimes\delta_y^* \sO_{\P^M}(\mu))=0.
$$
One has an exact sequence
\begin{multline*}
0\to \sO_{\tilde{\P}_y}(-(\mu+1)\cdot F)\otimes\delta_y^* \sO_{\P^M}(m)\\
\to \sO_{\tilde{\P}}(-\mu\cdot E)\otimes\delta^*
\sO_{\P}(m)|_{\tilde{\pi}^{-1}(y)}
\to \sO_{\P^M}(\mu) \to 0
\end{multline*}
and $H^1(\tilde{\pi}^{-1}(y), \sO_{\tilde{\P}}(-\mu \cdot
E)\otimes \delta^*\sO_{\P}(m)|_{\tilde{\pi}^{-1}(y)} )=0.$
By flat base change one finds
$$R^1\tilde{\pi}_*(\sO_{\tilde{\P}}(-\mu \cdot
E)\otimes \delta^*\sO_{\P}(m))=0.
$$
Moreover
$$\tilde{\pi}_*(\sO_{\tilde{\P}}(-\mu \cdot
E)\otimes \delta^*\sO_{\P}(m))|_y \to
\tilde{\pi}_*\sO_E(-\mu\cdot E)|_y \cong H^0(\P^M,\sO_{\P^M}(\mu))
$$
is an isomorphism. The inclusion
$$
S^{\mu} (\sF) \>>> \tilde{\pi}_* \delta^* \sO_{\P} (\mu)) \cong
S^\mu (\sE)
$$
factors through
$$
S^{\mu} (\sF) \> \subset >> \tilde{\pi}_* (\sO_{\tilde{\P}}
(-\mu \cdot E) \otimes \delta^* \sO_{\P} (\mu)).
$$
This map is an isomorphism. We know the surjectivity of
$$
\sF|_y \>>>
\tilde{\pi}_*(\sO_{\tilde{\P}}(-E)\otimes\delta^*\sO_{\P}(1))|_y
\cong H^0(\P^M,\sO_{\P^M}(1)),
$$
so for $\mu>1$ the morphism from $S^\mu(\sF)|_y$ to
$$
\tilde{\pi}_*(\sO_{\tilde{\P}}(-\mu\cdot E)\otimes\delta^*\sO_{\P}(\mu))|_y
\cong H^0(\P^M,\sO_{\P^M}(\mu))=S^\mu(H^0(\P^M,\sO_{\P^M}(1)))
$$
is surjective as well. By the choice of $s(B)$ one has
$$
\sO_{\P} (1) |_{s(B)} = \sO_{s(B)}
\mbox{ \ \ and \ \ }
\delta^* \sO_{\P} (1) |_E = \sO_E.
$$
Starting with $\mu = m$, assume by descending induction that
$$
H^0 (Y, \tilde{\pi}_* (\sO_{\tilde{\P}} (- \mu \cdot E) \otimes
\delta^* \sO_{\P} (m)) \otimes \sA^{-1} ) =0.
$$
Since
\begin{multline*}
H^0 (E, \tilde{\pi}_* (\sO_{\tilde{\P}} (- (\mu-1) \cdot E)
\otimes \delta^* \sO_{\P} (m)) \otimes \sA^{-1} |_E )\\
= H^0 (E, \tilde{\pi}_* (\sO_{\tilde{\P}} (- (\mu-1) \cdot E)
\otimes \delta^* \sO_{\P} (\mu -1)) \otimes \sA^{-1} |_E )\\
= H^0 (B, S^{\mu-1} (\sF) \otimes \sA^{-1} |_B )=0
\end{multline*}
one finds
$$
H^0 (Y, \tilde{\pi}_* (\sO_{\tilde{\P}} (- (\mu-1) \cdot E)
\otimes \delta^* \sO_{\P} (m)) \otimes \sA^{-1} =0 ,
$$
as well.
\end{proof}

\noindent
{\it Proof of \ref{7.1}.} \ Let us fix some $J \subseteq \{ 1,
\ldots , n\}$ with $S_J \neq \emptyset$. We will write $S_{\emptyset} =
Y$. The sheaf $S^m (T_Y (- {\rm log} \ D)) \otimes \sA$ is ample
with respect to $U_1 = Y - D$. Since $S + D$ is a normal
crossing divisor, $S_J \cap U_1 \neq \emptyset$ and since $\dim (S_J)
\geq 1$,
$$
H^0 (S_J, S^m (\Omega^{1}_{Y} ({\rm log} \ D)) \otimes \sA^{-1}
|_{S_J} ) =0.
$$
Assume, by induction on $\rho$, that
$$
H^0 (S_{J'} , S^m (\Omega^{1}_{Y} ({\rm log} (D+S_1 + \ldots
+S_{\rho-1} ))) \otimes \sA^{-1} |_{S_{J'}} ) =0,
$$
for all $m \geq 0$, and all $J' \subseteq \{ \rho , \ldots ,
\ell\}$ with $S_{J'} \neq \emptyset$. For $J \subseteq \{ \rho +1 ,
\ldots , \ell \}$ assume $T=S_{J} \neq \emptyset$. If $T_{\rho} =
S_{J \cup \{ \rho \}} = \emptyset$, i.e. if $S_{\rho} \cap T = \emptyset$,
then
$$
\Omega^{1}_{Y} ({\rm log} (D+S_1 + \ldots + S_{\rho} )) |_{S_J}
= \Omega^{1}_{Y} ({\rm log} (D+S_1 + \ldots +S_{\rho -1} ))
|_{S_J}
$$
and there is nothing to prove. Otherwise $T_{\rho} =
S_{\rho}|_T$ is a divisor and the restriction of
$$
0 \to \Omega^{1}_{Y} ({\rm log} (D+S_1 + \ldots +S_{\rho-1}))
\to \Omega^{1}_{Y} ({\rm log} (D+S_1 + \ldots +S_{\rho})) \to
\sO_{S_{\rho}} \to 0
$$
to $T$ remains exact. Hence for
$$
\sF = \Omega^{1}_{Y} ({\rm log} (D+S_1 + \ldots +S_{\rho-1}))
|_T \mbox{ \ \ and \ \ }
\sE = \Omega^{1}_{Y} ({\rm log} (D+S_1 + \ldots +S_{\rho}))
|_T
$$
$$
0 \to \sF \to \sE \to \sO_{T_S} \to 0
$$
is exact, $H^0 (T, S^m (\sF) \otimes \sA^{-1} |_T ) =0$ and
\begin{multline*}
H^0 (T_{\rho}, S^m (\sF) \otimes \sA^{-1} |_{T_S} )\\ = H^0 (S_{J
\cup \{ \rho\}} , S^m (\Omega^{1}_{Y} ({\rm log} (D+S_1 + \ldots
+ S_{\rho -1} )) \otimes \sA^{-1} |_{S_{J \cup \{ \rho\}}} ) =0.
\end{multline*}
Using \ref{7.5} we obtain
$$
H^0 (S_J, S^m (\Omega^{1}_{Y} ({\rm log} (D+S_1 + \ldots
+ S_{\rho } )) \otimes \sA^{-1} |_{S_J} ) =0.
$$ \qed
\begin{remark} \label{7.6}
The assumption ``$\sA$ ample'' was not really needed in the
proof of \ref{7.1}. It is sufficient to assume that
$$
\kappa (\sA |_{S_J} ) \geq 1, \ \mbox{for all} \ J \ \mbox{with}
\ S_J \neq \emptyset.
$$
\end{remark}

\bibliographystyle{plain}

\begin{thebibliography}{XXX} 
\bibitem{A-K} Abramovich, D., Karu, K.: Weak semi-stable
reduction in characteristic 0. Invent. math. {\bf 139} (2000)
241--273
\bibitem{B-V} Bedulev, E., Viehweg, E.: On the Shafarevich
conjecture for surfaces of general type over function fields.
Invent. Math. {\bf 139} (2000) 603--615
\bibitem{B-R} Br\"uckmann, P., Rackwitz, H.-G.: $T$-symmetrical
tensor forms on complete intersections. Math. Ann. {\bf 288} (1990)
627--635
\bibitem{Del} Deligne, P.: \'{E}quations diff\'{e}rentielles \`{a}
points singuliers r\'{e}guliers. Lecture Notes in Math. {\bf 163} (1970)
Springer, Berlin Heidelberg New York
\bibitem{E-V1} Esnault, H., Viehweg, E.:
Rev\^{e}tement cycliques. Algebraic Threefolds,
Proc. Varenna 1981. Springer Lect. Notes in Math. {\bf 947}
(1982) 241 - 250
\bibitem{E-V} Esnault, H., Viehweg, E.: Lectures on vanishing theorems.
DMV Seminar {\bf 20} (1992) Birkh\"auser, Basel Boston
\bibitem{Gri} Griffiths, P. (Editor): Topics in transcendental algebraic
geometry. Ann of Math. Stud. {\bf 106}, Princeton Univ. Press,
Princeton, NJ. (1984)
\bibitem{Gro} Grothendieck, A.: Techniques de construction et
th\'{e}or\`{e}mes d'existence en g\'{e}om\'{e}trie alg\'{e}brique, IV: Les
sch\'{e}mas de Hilbert. S\'{e}m. Bourbaki {\bf 221} (1960/61) In:
Fondements de la G\'{e}om\'{e}trie Alg\'{e}brique. S\'{e}m. Bourbaki,
Secr\'{e}tariat, Paris 1962
\bibitem{Gra} Grauert, H.: Mordells Vermutung \"uber rationale Punkte auf
algebraischen Kurven und Funktionenk\"orper. Publ. Math. IHES
{\bf 25} (1965), 131--149
\bibitem{Kaw} Kawamata, Y.: Minimal models and the Kodaira dimension
of algebraic fibre spaces. Journ. Reine Angew. Math. {\bf 363} (1985) 1--46
\bibitem{KKMS} Kempf, G., Knudsen, F., Mumford, D. and Saint-Donat, B.:
Toroidal embeddings I. Lecture Notes in Math. {\bf 339} (1973) Springer,
Berlin Heidelberg New York
\bibitem{Kov1} Kov\'acs, S.: Algebraic hyperbolicity of fine
moduli spaces. J. Alg. Geom. {\bf 9} (2000) 165--174
\bibitem{Kov2} Kov\'acs, S.: Logarithmic vanishing theorems
and Arakelov-Parshin boundedness for singular varieties.
preprint (AG/0003019), to appear in Comp. Math.
\bibitem{Kov3} Kov\'acs, S.: Families over a base with a
birationally nef tangent bundle. Math. Ann. {\bf 308} (1997) 347--359
\bibitem{Lu} Lu, S.S-Y.: On meromorphic maps into varieties of
log-general type. Proc. Symp. Amer. Math. Soc. {\bf 52} (1991) 305--333
\bibitem{Mig} Migliorini, L.: A smooth family of minimal
surfaces of general type over a curve of genus at most one is
trivial. J. Alg. Geom. {\bf 4} (1995) 353 - 361
\bibitem{O-V} Oguiso, K., Viehweg, E.: On the isotriviality of
families of elliptic surfaces. preprint (AG/9912100), to appear in J.
Alg. Geom.
\bibitem{Sch} Schmid, W.: Variation of Hodge structure: The singularities
of the period mapping. Invent. math. {\bf 22} (1973) 211--319
\bibitem{Sim} Simpson, C.: Higgs bundles and local systems. Publ. Math. I.H.E.S
{\bf 75} (1992) 5--95
\bibitem{Vie1} Viehweg, E.: Weak positivity and the additivity of the
Kodaira dimension for certain fibre spaces. In: Algebraic Varieties and
Analytic Varieties, Advanced Studies in Pure Math. {\bf 1} (1983) 329--353
\bibitem{Vie2}Viehweg, E.: Weak positivity and the additivity of the
Kodaira dimension for certain fibre spaces II. The local Torelli
map. In: Classification of Algebraic and Analytic
Manifolds, Progress in Math. {\bf 39} (1983) 567--589
\bibitem{Vie3} Viehweg, E.: Weak positivity and the stability of certain
Hilbert points. Invent. math. {\bf 96} (1989) 639--667
\bibitem{Vie} Viehweg, E.: Quasi-projective Moduli for
Polarized Manifolds. Ergebnisse der Mathematik, 3. Folge {\bf
30} (1995), Springer Verlag, Berlin-Heidelberg-New York
\bibitem{V-Z} Viehweg, E., Zuo K.: On the isotriviality of families
of projective manifolds over curves. preprint (AG/0002203). to appear in
J. Alg. Geom.
\bibitem{V-Z2} Viehweg, E., Zuo K.: On the Brody
hyperbolicity of moduli spaces for canonically
polarized manifolds. preprint (AG/0101004)
\bibitem{Zuo} Zuo, K.: On the negativity of kernels of
Kodaira-Spencer maps on Hodge bundles and applications.  Asian
J. of Math.  {\bf 4} (2000) 279--302
\end{thebibliography}

\end{document}